\documentclass[11pt]{article}
\input{pstricks.tex}
\input{mssymb.sty}

\newtheorem{thm}{Theorem}
\newtheorem{lem}{Lemma}
\newtheorem{cor} {Corollary}
\newtheorem{df} {Definition}

\newcommand {\simvol} {\parallel M,\partial M\parallel}
\newcommand {\hypvol} {\frac{1}{V_n}Vol\left(M\right)}
\newcommand {\chain} {\sum_{j=1}^ra_j\sigma_j}
\newcommand {\intschain} {\sum_{j=1}^ra_j\Phi_x\left(\sigma_j\right)}
\newcommand {\imagchain} {\cup_{j=0}^r im\left(\partial\sigma_j\right)}

\newcommand {\ch} {\overline{{\Bbb H}^n}}

\newcommand {\thick} {M_{\left[\epsilon,\infty\right]}}

\newcommand {\F} {\mathcal{F}}
\newcommand {\M} {\mathcal{M}}
\newcommand {\C} {\mathcal{C}}
\newcommand {\mH} {\mathcal{H}}
\newcommand {\sam} {\sum_{j=1}^ra_j}
\newcommand {\dura} {\cup_{j=0}^n A_j}
\newcommand {\durb} {\cup_{j=0}^n B_j}
\newcommand {\ver} {\cap_{j=0}^n\Gamma\backslash G - B_j}
\newcommand {\alv} {algvol\left(\sigma_j\right)}

\newenvironment{pf}{{\it Proof:}\quad}{\hfill Q.E.D. }

\begin{document}
\title{Efficient fundamental cycles of cusped hyperbolic manifolds }

\author{Thilo Kuessner\\
Mathematisches Institut, Universit\"at M\"unchen, 80333 M\"unchen, Germany}
\maketitle

\begin{abstract}
 
Let $M$ be a manifold (with boundary) of dimension $\ge3$, such that
its interior admits a hyperbolic metric of finite volume.
We discuss the possible
limits arising from sequences of relative fundamental
cycles approximating the simplicial volume $\parallel M,\partial M\parallel$, using ergodic theory of unipotent actions. 
As applications, we extend results of Jungreis and Calegari from closed
hyperbolic to finite-volume
hyperbolic manifolds:
 
a) strict subadditivity of simplicial volume with respect to isometric glueing
along geodesic surfaces, and
 
b) nontriviality of the foliated Gromov norm for "most"
foliations with two-sided branching.
\end{abstract}
\maketitle
\section{Introduction}
 
Gromov defined the simplicial volume $\parallel M,\partial M\parallel$
of a manifold $M$ as
the "minimal cardinality of a triangulation with real coefficients".
That means,
for an n-dimensional compact, connected, orientable
manifold $M$ with (possibly empty) boundary $\partial M$, define        
$$\simvol:=inf\left\{\sum_{i=1}^r\mid
a_i\mid:\sum_{i=1}^ra_i\sigma_i \mbox{ represents} \left[M,\partial
M\right]\right\}.$$ Here, $\left[M,\partial M\right]\in
H_n\left(M,\partial M;{\Bbb R}\right)$ is the image of a generator
of $H_n\left(M,\partial M;{\Bbb Z}\right)\simeq{\Bbb Z}$ under the canonical
homomorphism $H_n\left(M,\partial M;{\Bbb Z}\right)\rightarrow H_n\left(M,\partial M;
{\Bbb R}\right)$.
 
The simplicial volume quantifies the topological complexity of a
manifold. 
It is nontrivial 
if $int\left(M\right)$ admits a complete metric of sectional curvature $\le -a^2$ and finite volume.
In particular the Gromov-Thurston theorem (\cite{gro},\cite{thu})
states for finite-volume
hyperbolic manifolds
$\simvol=\hypvol$, where $V_n$ is
the volume of a regular ideal
simplex in ${\Bbb H}^n$.
This exhibits hyperbolic volume as a homotopy invariant, complementing the 
Chern-Gau\ss -Bonnet theorem, which implies homotopy invariance of hyperbolic volume
for even-dimensional manifolds. 
Homotopy invariance of hyperbolic volume was used by Gromov to give a more topological proof of Mostow's
rigidity theorem.
In the meantime, various more general rigidity theorems have been proved,
again using
the simplicial volume. 
 
On a finite-volume hyperbolic manifold, there does not exist a fundamental cycle actually having $l^1$-norm
$\hypvol$. However, there is a measure cycle,
supported on the set of regular ideal simplices, the so-called
smearing of a regular ideal simplex, having this norm:
after identifying the set of (ordered)
regular ideal simplices with $Isom\left({\Bbb H}^n\right)=                      
Isom^+\left({\Bbb H}^n\right)\cup Isom^-\left({\Bbb H}^n\right)$, it is the signed measure
$\frac{1}{2V_n}\left(Haar-r^*Haar\right)$, where $Haar$ is the Haar measure on $Isom^+\left({\Bbb H}^n\right)$
and $r$ is an orientation-reversing isometry.
This measure cycle can be approximated by
authentic
singular chains, i.e., finite linear combinations of (nonideal) simplices
(and this proves the Gromov-Thurston theorem, cf. \cite{bp} for details of the
proof in the case of closed manifolds).

Technically, the main part of this paper is devoted to the question of
to which extent this construction is unique, i.e., whether there
exist sequences of fundamental cycles with $l^1$-norms converging to
$\hypvol$ which do not approximate Gromov's
smearing construction.
 
For closed manifolds of dimension $\ge3$, it was shown in \cite{jun} by Jungreis
that {\em any} such sequence must converge to Gromov's smearing cycle.
In this paper we extend this rigidity results to hyperbolic manifolds of finite volume
which are either of dimension $\ge 4$ or which are of dimension $3$ and not Gieseking-like (see definition 4.4). 

Moreover, we obtain restrictions on sequences of fundamental cycles with $l^1$-norms converging to
$\hypvol$ on (possibly Gieseking-like) finite-volume hyperbolic 
manifolds of dimension $\ge3$, which allow to conclude: if $F$ is a closed geodesic
hypersurface, then the limits of such sequences 
do invoke simplices intersecting $F$ 'transversally' (see definition 5.1). This property can 
actually be restated as $\parallel M_F,\partial M_F\parallel > \simvol$,
where $M_F$ is obtained by cutting $M$ along $F$.\\
As applications, we extend results of Jungreis and Calegari to hyperbolic manifolds with cusps.\\

{\bf Glueing along boundaries.} Consider manifolds $M_1,M_2$, a
homeomorphism $f:A_1\rightarrow A_2$ between subsets
$A_i\subset\partial M_i$, and let $M=M_1\cup_fM_2$ be the
glued manifold. In general, it is hard to compare $\parallel M\parallel$
to $\parallel M_1\parallel+\parallel M_2\parallel$. One can prove
"$\ge$" if the $A_i$ are incompressible and amenable, and even "=" if,
in addition,
they are connected components of $\partial M_i$,
cf. \cite{gro} and \cite{kue}.\\
 
{\bf Theorem 6.3}:
{\em Let $n\ge3$ and let $M_1, M_2$ be compact n-manifolds with
boundaries $\partial M_i=\partial_0M_i\cup\partial_1M_i$, such that
$M_i-\partial_0 M_i$ admit incomplete hyperbolic metrics of finite volume
with totally geodesic boundaries $\partial_1M_i$.
If $\partial_1 M_i$ are not empty, $f:\partial_1 M_1
\rightarrow\partial_1 M_2$ is an isometry and $M=M_1\cup_f M_2$, then}
$$\parallel M,\partial M\parallel < \parallel M_1,\partial M_1
\parallel +\parallel
M_2,\partial M_2\parallel.$$
The same statement holds if one glues only along some connected components of $\partial_1 M_i$. One also has an analogous statement if two totally geodesic boundary components of
the same hyperbolic manifold are glued by an isometry.                      
 
One point of interest in theorem 6.3 is that it serves, in the case of 3-manifolds, as a
main step for a general glueing inequality. In 
\cite{kue}, we prove:\\

{\bf Theorem}:
{\em For a compact 3-manifold $M$
$\parallel DM\parallel <
2\parallel M,\partial M\parallel$ holds\\
if and only if $\parallel\partial M\parallel >0$,
i.e., if $\partial M$ consists not only of spheres and tori.}\\
 
(Here, $DM$ is
the manifold obtained by glueing two
differently oriented copies of $M$ via the identity of $\partial M$. Note that
$\parallel DM\parallel\le 2\parallel M,\partial M\parallel$
trivially holds.) This theorem may be seen as a generalisation of theorem 6.3, saying that any efficient fundamental cycle on 
a 3-manifold with ${\Bbb Z}_2$-symmetry has to intersect the fixed point set 'transversally'.
It is maybe worth pointing out that for the proof of this theorem in \cite{kue} we need to have theorem 6.3 also 
for the case of cusps. 
 
Another (direct) corollary from theorem
6.3 and Mostow rigidity is that (under the assumptions of theorem 6.3),
in dimensions $\ge4$, we get the same inequality for
{\em any} homeomorphisms $f$. This theorem seems to be hardly available
by topological methods.
The analogous statement in dimension 3
was recently shown to be wrong by Soma (\cite{soma}).
He proved: if $M_1,M_2$ are hyperbolic 3-manifolds of totally geodesic
boundary and $f:\partial M_1\rightarrow\partial M_2$ is pseudo-Anosov,
then $\lim_{n\rightarrow\infty}\parallel M_1\cup_{f^n}
M_2\parallel=\infty$.\\

{\bf Foliated Gromov norm.}
The Gromov norm of a foliation/lamination $\F$ on
a manifold $M$, as introduced in \cite{cal}, is
$$\parallel M,\partial M\parallel_{\F}:=inf\left\{\sum_{i=1}^r\mid
a_i\mid:\sum_{i=1}^ra_i\sigma_i \mbox{ represents } \left[M,\partial
M\right], \sigma_i\mbox { transverse to }\F\right\}.$$
The difference $\parallel M,\partial M\parallel_{\F}-\parallel M,\partial M\parallel$ 
seems to quantify the amount of branching of the leaf space.
Calegari proved:\\
- $\parallel M\parallel_{\F}=\parallel M\parallel$, when the leaf space is branched in
at most one direction, and\\
- $\parallel M\parallel_{\F}>\parallel M\parallel$
for asymptotically separated laminations
of closed hyperbolic manifolds of dimension $\ge3$.\\
The first statement generalizes easily to manifolds with boundary.
We extend the second statement as follows:\\
 
{\bf Theorem 7.5}:
{\em Assume that the interior of $M$ is a hyperbolic n-manifold of finite
volume. If $n\ge 3$ and $M$                                                        
is not Gieseking-like (definition 4.4), and
if $\F$ is an asymptotically
separated lamination, then
$$\parallel M,\partial M\parallel<\parallel M,\partial
M\parallel_{\F}.$$}\\
\\
We want to outline the content of this paper.
In chapter 3 we give a definition
of "efficient fundamental chains",
exhibit them as signed measures $\mu$ on the space of
regular ideal simplices,
show that they
are {\em absolute}
cycles (the boundary "escapes to infinity" for $\epsilon\rightarrow 0$),
and derive ergodic decompositions of
$\mu$ with respect to certain groups of reflections.
Such (different)
decompositions exist associated to all vertices of a fixed simplex $\Delta_0$. We
show that the ergodic decomposition corresponding to the $i$-th vertex of
$\Delta_0$ uses only the Haar measure and measures determined on the set
of simplices having $i$-th vertex in a parabolic fixed point of
$\Gamma$.
 
This is used in chapter 5 to prove theorem 5.3: if $F$ is a closed, totally geodesic
hypersurface in a finite-volume hyperbolic manifold of dimension $\ge3$ and
$\mu$ is an efficient fundamental cycle, then $\mu\left(S_F\right)>0$, where
$S_F$ is the set of simplices intersecting $F$ transversally.
To give a rough explanation of the proof:
the Haar measure does not vanish on $S_F$, hence $\mu\left(S_F\right)=0$
would imply that $\mu$ is determined on
the set of simplices with {\em all} vertices in
parabolic points, contradicting the fact that
it must invoke simplices with
faces in the cuspless hypersurface $F$.\\                                
In chapter 6, theorem 6.3 is derived from theorem 5.3. Chapter 7 is devoted to the foliated Gromov norm and the proof of theorem 7.5.
 
The simplicial volume of a nonorientable, disconnected manifold
is the sum over the connected components of half of the
simplicial volumina of the orientation coverings.
We will give all proofs for connected, oriented manifolds, since all statements generalise directly.
This includes that the orientations of glued manifolds
are understood
to fit together.
 
I would like to thank Michel Boileau and Bernhard Leeb for
discussions about the content of this paper.

\section{Preliminaries}                      
\subsection{Volume of straight simplices}
A simplex in hyperbolic space ${\Bbb H}^n$, with vertices $p_0,\ldots,p_i$, is called straight if it is the barycentric
parametrization of the geodesic simplex with vertices $p_0,\ldots, p_i$.

Given two regular ideal (straight) n-simplices $\Delta_0$ and $\Delta$ in ${\Bbb H}^n$,
with fixed orderings of their vertices, there is
a unique $g\in Isom\left({\Bbb H}^n\right)$ mapping $\Delta_0$ to $\Delta$. \\
Hence, fixing a reference simplex $\Delta_0$, we have an
$Isom\left({\Bbb H}^n\right)$-equivariant bijection between the set of ordered regular
ideal
n-simplices and
$Isom\left({\Bbb H}^n\right)$, this bijection being unique up to the choice of $\Delta_0$, i.e., up to multiplication with a fixed element of $Isom\left({\Bbb H}^n\right)$. 

As another consequence, all regular ideal n-simplices in ${\Bbb H}^n$ have the
same volume, to be denoted $V_n$.

By
\cite{hamu}, any straight n-simplex $\sigma$ in ${\Bbb H}^n$ satisfies
$Vol\left(\sigma\right)\le V_n$ and equality is achieved only for
regular ideal simplices.

\subsection{Ergodic decomposition}
For a topological space $X$, we consider Radon measures $\mu$ on $X$. This are, by definition, elements of $C_c^*\left(X\right)$, the dual 
of the space of compactly supported continuous functions. They have a decomposition $\mu=\mu^+-\mu^-$ with $\mu^+,\mu^-$ nonnegative Radon measures. (We will refer to $\mu$ as signed measure and to $\mu^\pm$ as measures.) A probability measure on $X$
is a measure $\mu$ with $\mu\left(X\right)=1$.\\
Let a group $G$ act on a topological space $X$. A probability measure $\mu$ is called ergodic if any $G$-invariant set has measure $0$ or $1$.
Denote by $\mathcal{E}$ the set of ergodic G-invariant probability measures
on $X$. \\
Let $\mathcal{A}$ be the weak measure class induced by the measure class on $X$, i.e., the
smallest $\sigma$-algebra $\mathcal{A}$ on $\mathcal{E}$ 
such that 
for all Borel sets $A\subset X$ the application $f_A:\mathcal{E}\rightarrow {\Bbb R}$
defined by $$f_A\left(\mu\right):=\mu\left(A\right)$$
is measurable.
\begin{lem}: Let a group G act on a complete separable metric space X. If there
exists a
G-invariant probability measure on $X$,
then the set $\mathcal{E}$ of ergodic G-invariant measures
on $X$ is not empty and there is a decomposition map
$\beta:X\rightarrow\mathcal{E}$.\end{lem}
Here, a decomposition map is a $G$-invariant map $\beta:X\rightarrow\mathcal{E}$, 
which is\\
- measurable with respect to $\mathcal{A}$,\\ 
- satisfies
$e\left(\left\{x\in X:\beta\left(x\right)=e\right\}\right)=1$ for all $e\in
\mathcal{E}$, and\\
- for all $G$-invariant probability measures
$\mu$ and Borel sets $A\subset X$ the following equality holds:
$$\mu\left(A\right)=
\int_X\beta\left(x\right)\left(A\right)d\mu\left(x\right).$$
For a proof of lemma 2.1, see theorem 4.2 in \cite{var}.\\          
\\
For later reference we state the following lemma, part (i) of which is known as
Alaoglu's
theorem, whereas a proof of part (ii) can be found in lemma 3.2 of \cite{dan2}.
\begin{lem}: (i) Any weak-*-bounded sequence of signed Radon measures on
a locally compact metric space has an accumulation point in the weak-*-topology.\\
(ii) If $\mu$ is the weak-*-limit of a sequence $\mu_n$ of measures on a space $X$,      and $U\subset X$ is an open subset, then
$\mu\left(U\right)\le\liminf\mu_n\left(U\right)$. \end{lem}  
Moreover, we recall that the support of a measure $\mu$ on $X$ is defined as the
complement of the largest open set $U\subset X$ with the property $\mu\left(U\right)=0$.
\subsection{Measure homology}
The following explanations are not necessary (from a logical point of view) for our arguments, but may be helpful to 
understand the framework.
For a topological space $X$, let $C^0\left(\Delta^k,X\right)$
be the space of singular $k$-simplices in $X$, topologized by the
compact-open-topology. For a signed measure $\mu$ on $C^0\left(\Delta^k,X\right)$,
one has its decomposition $\mu=\mu^+-\mu^-$ as difference of two (non-negative) Borel measures, and one defines its total variation
as $\parallel \mu\parallel=\int d\mu^+ +\int d\mu^-$.

Let ${\C}_k\left(X\right)$ be the vector space of all
signed measures $\mu$ on $C^0\left(\Delta^k,X\right)$ which have compact support and finite
total variation. (We assume finite total variation because we want $\parallel .\parallel$ to define a norm on ${\C}_k\left(X\right)$.The condition 'compact support' is 
imposed because otherwise the map $j_*:H_*\left(X;{\Bbb R}\right)\rightarrow 
{\mathcal{H}}_*\left(X\right)$ defined below would, in general, not be surjective, see \cite{thu}, 6.1.7 for examples of this phenomenon.)
Let $\eta_i:\Delta^k\rightarrow\Delta^{k-1}$
be the i-th face map. It induces a map $\partial_i=\left(\eta_i^*\right)_*:
{\C}_k\left(X\right)
\rightarrow {\C}_{k-1}\left(X\right)$. We define the boundary operator $\partial:=\sum_{i=0}^k\partial_i$, to make ${\C}_*\left(X\right)$ a chain complex. We
denote the homology groups of this chain complex by ${\mH}_*\left(X\right)$.

We have an obvious inclusion $j: C_*\left(X\right)\rightarrow
{\C}_*\left(X\right)$, where 
$ C_*\left(X\right)$ are the singular chains, considered as
finite linear combination of atomic measures. Clearly, $j$ is a chain map.
Zastrow's theorem 3.4. in \cite{zas} says that we get an isomorphism 
$j_*:H_*\left(M\right)\rightarrow {\mH}_*\left(M\right)$ if $M$ is a smooth manifold (but not for arbitrary topological spaces $X$).

The $l^1$-norm on $C_*\left(M\right)$ extends to the total variation $\parallel.\parallel$ on
${\C}_*\left(M\right)$, and we get an induced pseudonorm on 
${\mathcal{H}}_*\left(M\right)$. Thurston conjectured in \cite{thua} that the isomorphism
$j_*$ should be an isometry for this pseudonorm. There seems not to exist a proof of this general conjecture
so far, but if $M$ is a closed hyperbolic n-manifold,
it follows easily
from the identity $\parallel M\parallel=\frac{1}{V_n}Vol\left(M\right)$ (\cite{bp},\cite{gro})
that $j_n:H_n\left(M\right)\rightarrow {\mathcal{H}}_n\left(M\right)$ is
an isometry.\\

\subsection{Intersection numbers}
\begin{df}: Let $M$ be an oriented differentiable n-manifold.
For an immersed differentiable
n-simplex $\sigma:\Delta^n\rightarrow M$, and $x\in M$,
define
$$\Phi_x\left(\sigma\right)=
\sum_{y\in\sigma^{-1}\left(x\right)}
\mbox{sign }\det d\sigma\left(y\right).$$ 
For a singular chain
$c=\chain$, let
$\Phi_x\left(c\right)=
\intschain$. \end{df}  
\begin{lem}: Let $M$ be a connected, oriented, smooth, noncompact n-manifold, $M^\prime$ an  
n-submanifold with boundary, such that $M - M^\prime$ is compact.
Let $c=\chain$ be a smooth singular n-chain
representing the relative fundamental class
$\left[M, M^\prime\right]$. Assume that all $\sigma_j$ are immersed smooth
n-simplices. Then $\Phi_x
\left(c\right)= 1$ holds for almost all $x\in M - M^\prime$. \end{lem}
\begin{pf} 
Let $K=\imagchain$. $K$ is
of measure zero, by Sard's lemma.

We want to show that $\Phi_x\left(
c\right)$, as a function of $x$, is constant on 
$M - \left(M^\prime\cup K\right)$.
It is obvious that it is locally constant
on $M - \left(M^\prime\cup K\right)$, since all $\sigma_i$ are immersed.
It remains to prove: for all $x\in K\cap\mbox{int}\left(
M - M^\prime\right)$,
there is a neighborhood $U$ of $x$ in $M$ such that
$\Phi_.\left(c\right)$ is constant on $U\cap \left(M - K\right)$.
 
The point $x$ is contained in the image of finitely many (n-1)-simplices $\kappa_1,
\ldots,\kappa_k$,
which are boundary
faces of some $\sigma_{i_1},\ldots,\sigma_{i_k}$.
(Note that the $\sigma_{i_j}$'s
need not be distinct and that there might be further $\sigma_i$'s
containing $x$ in the interior of their image.)
Since $\partial\chain$ invokes only simplices whose image
is contained in $M - M^\prime$, we necessarily have that all $\sigma_{i_1},
\ldots,\sigma_{i_k}$ cancel each other, i.e., there is a partition of
$\left\{i_1,\ldots,i_k\right\}$ in some subsets,
such that for each of these subsets of indices
the sum of the corresponding coefficients $a_{i_j}$,
multiplied with a sign according to
orientation of $\sigma_{i_j}$,
adds up to zero.
This implies that $\Phi_.$ is constant in the intersection of
a small neighborhood of $x$ with the complement of $K$ and, hence, also
constant on all of $M - \left(M^\prime\cup K\right)$.                                                  
 
We now prove that this constant does
not depend on the representative of
the relative fundamental class. This implies that the constant must be 1,
since one can choose a triangulation as representative of the
relative fundamental class. 
 
If $c$ and $c^\prime$
are different representatives
of $\left[M,M^\prime\right]$, we have that $c-c^\prime=\partial w+t$
for some $w\in C_{n+1}\left(M\right)$ and $t\in C_n\left(
M^\prime\right)$. 
Because $\partial w$ is a cycle, the same
argument as above gives that $\Phi_.\left(\partial w
\right)$
is a.e.\ constant on {\em all} of $M$. This constant must be zero, since
$\partial w$
has compact support in the noncompact manifold
$M$. That means that
$\Phi_x\left(c\right)-\Phi_x\left(c^\prime
\right)=\Phi_x\left( t\right)$
for almost all $x\in M$. 
But $\Phi_x\left(t\right)=0$ for all
$x\in int\left(M - M^\prime\right)$.\end{pf}\\       

\subsection{Convergence of fundamental cycles }

A major point of the next chapter will be to consider limiting objects of sequences of
relative fundamental cycles  
of a finite-volume hyperbolic manifold $M$ with $l^1$-norms converging to the simplicial 
volume. Since straight simplices have volume smaller than $V_n$, there do not exist relative fundamental cycles
actually having $l^1$-norm equal to $\hypvol$.
Hence, the limits of such sequences can not be just singular chains. 
What we are going to do is to embed the singular chain complex into a larger
space, where any bounded sequence has accumulation points. A straightforward idea
would be to use the inclusion $j:C_n\left(M\right)\rightarrow {\C}_n\left(M\right)$
and to consider weak-* accumulation points in ${\C}_n\left(M\right)$.
This works perfectly well, however it is easy to see that the weak-* limits
are just trivial measures. The reason is roughly 
the following: a singular chain with $l^1$-norm close to $\hypvol$
has to have a very large part of its mass on simplices $\sigma$ with $vol\left(str\left(\sigma\right)\right)$ quite close to $V_n$. If we consider a compact set of simplices, it will have some upper bound
(better than $V_n$) on $vol\left(str\left(.\right)\right)$. Hence, it will contribute 
very little to an almost efficient fundamental cycle, and the limiting measure will actually 
vanish on this set of simplices. 

Therefore, to get nontrivial accumulation points, we are obliged to consider the larger
space of simplices which might be ideal, i.e., whose lifts to ${\Bbb H}^n$ might have vertices 
in $\partial_\infty {\Bbb H}^n$. This, however, raises another problem: the space of ideal
simplices in $M=\Gamma\backslash {\Bbb H}^n$ is not Hausdorff, and there is no theorem 
guaranteeing existence of weak-* accumulation points for signed measures on non-Hausdorff
spaces.
 
(If $M=\Gamma\backslash{\Bbb H}^n$ has finite volume, then the action of $\Gamma$ on 
$\partial_\infty{\Bbb H}^n$ has dense orbits. Thus the quotient $\Gamma\backslash SS_i\left(M\right)$ can 
not be Hausdorff as long as $SS_i\left(M\right)$ contains degenerate simplices. We will show
in section 2.5.3, however, that the
action of $\Gamma$ on the subspace of nondegenerate $n$-simplices is properly discontinuous, i.e., after throwing away the
degenerate simplices we get a Hausdorff quotient. A similar idea seems to have been exploited in the proof of lemma 2.2 in \cite{jun} where the author
restricted to a compact subset of $SS_n\left(M\right)$, i.e., to simplices with a lower volume bound.)\\

\subsubsection{Straightening alternating chains}

The symmetric group $S_{n+1}$ acts on the standard
n-simplex $\Delta^n$: any permutation
$\pi$ of vertices can be realised by a 
unique affine map $f_\pi:\Delta^n\rightarrow\Delta^n$. For a singular simplex $\sigma:\Delta^n\rightarrow M$ let
$alt\left(\sigma\right):=\sum_{\pi\in S_{n+1}}sgn\left(\pi\right)
\sigma f_{\pi}$, and for a singular chain $c=\chain$ define
$alt\left(c\right):=\frac{1}{\left(n+1\right)!}\sum_{j=1}^r a_j alt\left(\sigma_j\right)$. Clearly,
$\parallel alt\left(c\right)\parallel\le\parallel c\parallel$.                               
 
For a simplex $\sigma$ in ${\Bbb H}^n$, we denote
by Str($\sigma$) the straight simplex with the same vertices as
$\sigma$ (as in section 2.1.).
A straight simplex in a hyperbolic manifold $M=\Gamma\backslash {\Bbb H}^n$ is
the image of a straight simplex in ${\Bbb H}^n$ under the projection
$p:{\Bbb H}^n\rightarrow \Gamma\backslash {\Bbb H}^n=M$. For a simplex $\sigma$ in
$M$, its
straightening $Str\left(\sigma\right)$ is defined as
$p\left(Str\left(\tilde{\sigma}\right)\right)$, where $\tilde{\sigma}$
is a simplex in ${\Bbb H}^n$ projecting to $\sigma$.
Since straightening in ${\Bbb H}^n$ commutes with isometries, the definition of
$Str\left(\sigma\right)$ does not depend on the choice of $\tilde{\sigma}$.
 
Finally, the straightening of a singular chain c=$\chain$ is defined
as \\
$Str\left(c\right)=\sum_{j=1}^ra_j Str\left(\sigma_j\right)$. $Str\left(c\right)$ is homologous to $c$, and
clearly \\
$\parallel Str\left(c\right)\parallel
\le\parallel c\parallel$ for 
any $c\in C_*\left(M\right)$. ($Str\left(c\right)$ may possibly have smaller
norm than $c$, since different simplices can have the same straightenings.)

If $M^\prime\subset M$ is a convex subset (meaning that $\sigma\subset M^\prime$ implies $str\left(\sigma\right)\subset M^\prime$), then $Str:C_*\left(M,M^\prime\right)\rightarrow C_*\left(M,M^\prime\right)$ is well-defined.\\

\subsubsection{Nondegenerate chains}

Let $M$ be a hyperbolic manifold.
We call a {\em straight} i-simplex $\sigma:\Delta^i\rightarrow N$ degenerate if two of its vertices are 
mapped to the same point, nondegenerate otherwise. \\
For $M^\prime\subset M$ a convex subset of $M$,
we consider $algvol:C_n\left(M,M^\prime\right)\rightarrow R$ which maps $\sigma\in C_n\left(M\right)$
to the algebraic volume (see \cite{bp}, p.107) of $str\left(\sigma\right)\cap\left(M-M^\prime\right)$. (Since $M^\prime$ is 
convex, $algvol$ is well-defined on the relative chain complex.) It follows from Stokes theorem that we get an induced map
$algvol_*:H_n\left(M,M^\prime;{\Bbb R}\right)\rightarrow {\Bbb R}$.

\begin{lem}: Let $M$ be a hyperbolic n-manifold, $M^\prime$ a convex subset
such that $algvol:H_n\left(M,M^\prime;{\Bbb R}\right)\rightarrow {\Bbb R}$ is an isomorphism.
Let
$\sum_{i\in I}a_i\sigma_i\in C_n\left(M,M^\prime;{\Bbb R}\right)$ be a straight relative n-cycle.
Then there is a subset of indices $J\subset I$ such that all $\sigma_j$
with $j\in J$ are non-degenerate and 
$\sum_{j\in J}a_j\sigma_j$ is relatively homologous
to $\sum_{i\in I}a_i\sigma_i$.\end{lem}
\begin{pf}
Let $K:=\left\{k\in I: \sigma_k\mbox{ degenerate }\right\}$ be the set of indices of degenerate simplices occuring in $\sum_{i\in I}a_i\sigma_i$.
We claim that
$\sum_{k\in K} a_k\sigma_k$ is a relative cycle.
Indeed, the degenerate faces of $\sum_{k\in K}a_k\sigma_k$ cancel each other
(relatively),
since they cancel in $\partial\left(\sum_{i\in I}a_i\sigma_i\right)$
and they can not cancel against faces of nondegenerate simplices.
Moreover, 
the nondegenerate faces of degenerate simplices cancel anyway: if
$\left(a,v_1,\ldots,v_n\right)$ and $\left(b,v_1,\ldots,v_n\right)$ are
nondegenerate faces of a degenerate simplex, then necessarily $a=b$.
Thus this face contributes twice to the boundary, with opposite signs. 

We have obtained that $\sum_{k\in K}a_k\sigma_k$ is a relative cycle. 
But, since all $\sigma_k$ are degenerate, they have vanishing volume, and
we have that the relative 
homology class $\left[\sum_{k\in K}a_k\sigma_k\right]\in ker\left(algvol_*\right)=0$ (since $algvol_*$ is an isomorphism, by assumption), i.e., $\sum_{k\in K}a_k\sigma_k\in ker\left(algvol\right)=0$ is a relative boundary. Then choose $J=I-K$.
\end{pf}\\

In conclusion, if $M^\prime\subset M$ convex and $n=dim\left(M\right)$,
then to any relative n-cycle $c\in C_n\left(M,M^\prime; {\Bbb R}\right)$ 
we find 
$c^\prime\in C_n\left(M,M^\prime;{\Bbb R}\right)$ homologous to $c$ in $C_*\left(M,M^\prime;{\Bbb R}\right)$, such that 
$\parallel c^\prime\parallel \le\parallel c\parallel$ and $c^\prime$ is an alternating linear combination of nondegenerate straight simplices.\\
 
\subsubsection{Straight chains as measures}
 
We explained in 2.3 that singular chains may be considered as measures
on the space of singular simplices, thus getting a homomorphism 
$C_*\left(M;{\Bbb R}\right)\rightarrow {\C}_*\left(M\right)$. 
As we said, to get nontrivial results, we should consider not only
${\C}_*\left(M\right)$, but measures on the space of possibly ideal
simplices. Since it is hard
to prove existence of accumulation points in this measure space,
we will consider measures on smaller sets of simplices.

Let $M$ be a hyperbolic manifold.
The set of nondegenerate, possibly ideal,
straight i-simplices in $M=\Gamma\backslash {\Bbb H}^n$ is
$$SS_i\left(M\right):=\Gamma\backslash \left\{\left(p_0,\ldots,p_i\right):
p_0,\ldots,p_i\in\ch, p_j\not=p_k\mbox{ if }j\not=k\right\},$$
where $g\in\Gamma$ acts by
$g\left(p_0,\ldots,p_n\right)=\left(gp_0,\ldots,gp_n\right)$.
 
Denote ${\M}\left(SS_i\left(M\right)\right)$ the space of signed regular measures on $SS_i\left(M\right)$. Straight 
singular chains $c=\chain\in C_i\left(M;{\Bbb R}\right)$, with all $\sigma_j$ nondegenerate,
can be considered as discrete signed measures
on $SS_i\left(M\right)$
defined by
$$c\left(B\right)=\sum_{\left\{j:\sigma_j\subset
B\right\}}\mid a_j\mid$$ for any Borel set $B\subset SS_i\left(M\right)$.

Let $n=dim\left(M\right)$.
To apply Alaoglu's theorem to ${\M}\left(SS_n\left(M\right)\right)$,
we need to know that $SS_n\left(M\right)$ is locally compact (which is obvious)
and metrizable.

\begin{lem}:
Let $M$ be a hyperbolic manifold of dimension $n\ge3$. Then $SS_n\left(M\right)$
is metrizable.\end{lem}
\begin{pf}
We have to show
that $\Gamma$-orbits on $\Pi_{j=0}^n\ch -D$ are closed, $D$ being the set of
degenerate straight simplices.
On the complement of $\Pi_{j=0}^n\partial_\infty\ch$ this follows from
proper discontinuity of the $\Gamma$-action on ${\Bbb H}^n$.
 
Now we assume $n\ge 3$.
To any n-tuple $\left(v_0,\ldots,v_{n-1}\right)
\in\Pi_{j=0}^{n-1}\partial_\infty\ch$
of {\em distinct} points
corresponds a unique $v_n\in\partial_\infty\ch$ such that
$\left(v_0,\ldots,v_n\right)$ is a positively oriented regular ideal n-simplex. (If $n=2$, then $v_n$ is not uniquely determined.)
Together with 2.1, we get
a $\Gamma$-equivariant
homeomorphism\\
$\Pi_{j=0}^{n-1}\partial_\infty
\ch - D\rightarrow Isom^+\left({\Bbb H}^n\right)$. Since $\Gamma\backslash{\Bbb H}^n$ is a manifold, we know that
$\Gamma$ acts properly discontinuously
on $Isom^+\left({\Bbb H}^n\right)$, thus also on $\Pi_{j=0}^{n-1}\partial_\infty
\ch - D$. This implies of course that it acts properly discontinuously on $\Pi_{j=0}^n\partial_\infty\ch - D$. Thus, $\Gamma$-orbits are closed.\end{pf}

\section{Degeneration}
\subsection{Efficient fundamental cycles}
For a closed hyperbolic
manifold $M$, we know that $\parallel M\parallel=\hypvol$. This means that, for any
$\epsilon >0$, there is some fundamental cycle $c_\epsilon$
satisfying
$\parallel c_\epsilon\parallel\le\parallel M\parallel +\frac{\epsilon}{V_n}$.
By 2.5.1 and 2.5.2, we can choose $c_\epsilon$ to be an alternating chain 
consisting of
nondegenerate straight simplices, without increasing the $l^1$-norm.
To speak about limits of sequences
of $c_\epsilon$, one has to regard them as elements of
some locally compact space,
namely the space of signed Radon measures on $SS_n\left(M\right)=
\Gamma\backslash\left(\Pi_{j=0}^n\ch - D\right)$ with the weak-*-topology,
as in 2.5.3.\\
 
Jungreis' results from \cite{jun}, for closed hyperbolic manifolds of dimension
$\ge3$, can be rephrased as follows:\\
- any sequence of $c_\epsilon$ as above,
with $\epsilon\rightarrow 0$,
converges,\\
- the limit is a signed measure $\mu$, which is supported
on the set of regular ideal simplices (to be identified with $Isom\left({\Bbb H}^n
\right)$), and\\
- up to a multiplicative factor one has
$\mu=\mu^+-\mu^-$ with $\mu^+$ the Haar measure on $Isom^+\left({\Bbb H}^n\right)$
and $\mu^-=r^*\mu^+$ for an arbitrary orientation
reversing $r\in Isom\left({\Bbb H}^n\right)$.\\
 
The aim of this chapter is to generalize these results to finite-volume hyperbolic manifolds.
For these cusped hyperbolic manifolds, there arises a technical problem:
we wish to consider chains representing the relative fundamental class of a
manifold with boundary, but we have a hyperbolic metric (and a notion of
straightening) only on the interior. In the following, we will get around
this problem and analyse the possible limits.\\

Let $M$ be a compact $n$-manifold with boundary $\partial M$ such that $int\left(M\right)$ carries a hyperbolic metric of finite volume.
With respect to this hyperbolic metric, denote $M_{\left[a,b\right]}:=\left\{x\in int\left(M\right): a\le inj\left(x\right)\le b\right\}$. 
It is a well-known consequence of the Margulis lemma (\cite{bp}, D.3.12.) that, for sufficiently small $\epsilon >0$, the '$\epsilon$-thin part' $M_{\left[0,\epsilon\right]}$ is a product neighborhood of $\partial M$, i.e., homeomorphic to $\partial M\times\left[0,\infty\right)$. Thus, one has a retraction $r_\epsilon$ from $M$ to the '$\epsilon$-thick part' $M_{\left[\epsilon,\infty\right]}$ which induces a homeomorphism of pairs $r_\epsilon:\left(M,\partial M\right)\rightarrow \left(\thick,\partial\thick\right)$ and, thus, an isomorphism $$r_{\epsilon *}:H_*\left(M,\partial M\right)\rightarrow H_*\left(\thick,\partial\thick\right).$$ (This applies to all $\epsilon<\epsilon_0$, where $\epsilon_0$ depends on $M$.)

It should be noted that $M_{\left[0,\epsilon\right]}$ is convex and that one has the isomorphism
$$ algvol:H_n\left(M,M_{\left[0,\epsilon\right]};{\Bbb R}\right)\rightarrow{\Bbb R}.$$
Convexity of $M_{\left[0,\epsilon\right]}$ implies that the
straightening homomorphism $$Str:C_*\left(int\left(M\right),M_{\left[0,\epsilon\right]}\right)\rightarrow
C_*\left(int\left(M\right),M_{\left[0,
\epsilon\right]}\right)$$ is well-defined and induces an isomorphism in relative homology. Moreover, there is the inclusion 
$$exc:C_*\left(\thick,\partial\thick\right)\rightarrow C_*\left(int\left(M\right),M_{\left[0,\epsilon\right]}\right),$$
which induces an isomorphism in homology by the excision theorem.
In conclusion, $$Str\left(exc\left(r_{\epsilon*} .\right)\right):C_n\left(M,\partial M;{\Bbb R}\right)\rightarrow C_n\left(int\left(M\right),M_{\left[0,\epsilon\right]};{\Bbb R}\right)$$ induces an isomorphism in homology and does not increase $l^1$-norms.\\

Let, for $\epsilon < \epsilon_0$, $c_\epsilon\in C_n\left(M,\partial M;{\Bbb R}\right)$ be some relative
fundamental cycle satisfying
$$\parallel c_\epsilon\parallel\le
\simvol +\frac{\epsilon}{V_n}. $$
By the above arguments, we may replace $c_\epsilon$ by \\
$Str\left(exc\left(r_{\epsilon *} c_\epsilon\right)\right)\in C_n\left(int\left(M\right),M_{\left[0,\epsilon\right]};{\Bbb R}\right)$ without increasing the $l^1$-norm. Abusing notation, we will continue to denote this new relative cycle by $c_\epsilon$.

\begin{df}: A signed measure $\mu$ on $SS_n\left(M\right)$
is called an efficient fundamental chain if there exists a sequence of $\epsilon$ with $\epsilon\rightarrow 0$ and a sequence of $c_\epsilon\in C_n\left(M,\partial M;{\Bbb R}\right)$ representing the relative fundamental class $\left[M,\partial M\right]$,  
which are alternating chains invoking only nondegenerate simplices and
which satisfy
$\parallel c_\epsilon\parallel\le
\simvol +\frac{\epsilon}{V_n}$, such that the sequence $Str\left(exc\left(r_{\epsilon*}c_\epsilon\right)\right)\in C_n\left(
int\left(int\left(M\right)\right),M_{\left[0,\epsilon\right]};{\Bbb R}\right)$ converges to $\mu$
in the weak-*-topology of ${\mathcal{M}}\left(SS_n\left(M\right)\right)$, the space of signed measures on the space of 
straight nondegenerate simplices.\end{df}
\begin{lem}: Assume that $M$ is a manifold of dimension $n\ge 3$, such that $int\left(M\right)$ admits a hyperbolic metric of finite volume.
Then there is at least one efficient fundamental chain.\end{lem}
\begin{pf} Considering some sequence of $c_\epsilon$ with
$\epsilon\rightarrow 0$, we may by 2.5 assume that the support of the $c_\epsilon$ consists of
only straight nondegenerate simplices.
$Str\left(exc\left(r_{\epsilon*}c_\epsilon\right)\right)$
may be regarded as a sequence of
signed measures on the locally compact metric space $SS_n\left(M\right)$,
see 2.5.3.
The sequence $c_\epsilon$
is bounded by its definition and, hence, lemma 2.2 and 2.6 guarantee the existence of
a weak-*-accumulation point $\mu$. (The condition $n\ge 3$ is needed to apply lemma 2.6.)
\end{pf}\\
\\

We recall that excision and straightening, as well as the homeomorphism
$r_\epsilon$ induce isomorphisms in relative homology. Hence,
any new $c_\epsilon$ represents the relative fundamental class in 
$H_n\left(int\left(M\right),M_{\left[0,\epsilon\right]};{\Bbb R}\right)$. 
As a special case of lemma 2.4 we have:
\begin{lem}: Let $c_\epsilon$
be a representative of the relative fundamental class
$\left[int\left(M\right),M_{\left[0,\epsilon\right]}\right]$. Then $\Phi_x
\left(c_\epsilon\right)= 1$ holds for almost all $x\in\thick$.\end{lem}
\begin{df} For a hyperbolic manifold $M$, let $S_\delta\subset SS_n\left(M\right)$ 
be the set of nondegenerate straight simplices $\sigma \in M$ with $vol\left(\sigma\right) < V_n-\delta$.\end{df}
\begin{lem}: An efficient fundamental chain $\mu$ is supported on \\
$SS_n\left(M\right) - S_0$, i.e., on the set
of straight simplices of volume $V_n$.\end{lem}
\begin{pf}
It suffices to show that $\mu\left(S_{\delta}\right)=0$ holds for
any $\delta > 0$. By lemma 2.2, (ii), and openness of $S_{\delta}$, this follows if we can prove $\lim_{\epsilon\rightarrow 0}c_\epsilon\left(S_{\delta}\right)=0$ for
any $\delta > 0$. Here, $c_\epsilon=\chain$ is the sequence from definition 3.1.

From lemma 3.3, we conclude
$\int_M\Phi_x\left(c_\epsilon\right)dvol\left(x\right)\ge
Vol\left(\thick\right)$.\\
But $\int_M\Phi_x\left(c_\epsilon\right)dvol\left(x\right)=
\sam \int_M\Phi_x\left(\sigma_j\right)dvol\left(x\right)=\sam\ \alv$,
where $algvol\left(\sigma_j\right)$ is
$Vol\left(\sigma_j\right)$ with
a sign according to orientation. As a consequence:
$$\sum\mid a_j\mid Vol\left(\sigma_j\right)\ge
Vol\left(\thick\right).$$
On the other hand, we want $c_\epsilon=\chain$ to satisfy
$V_n\sum \mid a_j\mid \le Vol\left(M\right)+\epsilon$.\\
Subtracting the two inequalities                                                    yields
$$\sum\mid a_j\mid \left(V_n-Vol\left(\sigma_j\right)\right)\le\epsilon+
Vol\left(M_{\left[0,\epsilon\right]}\right).$$
 
We get \\
$\hspace*{0.3in}\epsilon+Vol\left(M_{\left[0,\epsilon\right]}\right)\ge 
\sum\mid a_j\mid\left(V_n-Vol\left(\sigma_j\right)\right) \\
=\sum_{j:Vol\left(\sigma_j\right) \ge V_n-\delta}
\mid a_j\mid\left(V_n-Vol\left(\sigma_j\right)\right)+
\sum_{j:Vol\left(\sigma_j\right) < V_n-\delta}
\mid a_j\mid\left(V_n-Vol\left(\sigma_j\right)\right)\\
 \ge \sum_{j:Vol\left(\sigma_j\right) < V_n-\delta}
\mid a_j\mid\left(V_n-Vol\left(\sigma_j\right)\right)\\
 \ge \delta\sum_{j:Vol\left(\sigma_j\right) <
V_n-\delta}\mid a_j\mid\\
\hspace*{0.3in}  = \delta c_\epsilon\left(S_{\delta}\right)$.\\
Since
$\lim_{\epsilon\rightarrow 0}Vol
\left(M_{\left[0,\epsilon\right]}\right)=0$, we conclude $\lim_{\epsilon\rightarrow 0}c_\epsilon\left(S_{\delta}\right)=0$.\end{pf}\\
\\
We have just proved that, if $c_\epsilon$ is a representative of 
$\left[int\left(M\right),M_{\left[0,\epsilon\right]}\right]$
satisfying $V_n\parallel c_\epsilon\parallel\le Vol\left(M\right)+\epsilon$,
then $\lim_{\epsilon\rightarrow 0}c_\epsilon\left(S_\delta\right)=0$ holds for any $\delta>0$. Since $\parallel 
c_\epsilon^\pm\parallel\le\parallel c_\epsilon\parallel$ and $\int_M\Phi_x\left(c_\epsilon^\pm\right)\ge Vol\left(\thick\right)$, we can use the same argument to prove that 
$\lim_{\epsilon\rightarrow 0}c^\pm_\epsilon\left(S_\delta\right)=0$ holds for any $\delta>0$. This fact will be used in the proof of the following lemma.

\begin{lem}: Let $\mu$ be an efficient fundamental chain. Then $\mu\not=0$.\end{lem}
\begin{pf} Choose a continuous $f:SS_n\left(M\right)\rightarrow\left[0,1\right]$, which vanishes on some $S_{\delta^\prime}$ and
is $1$ on the complement of some $S_\delta$.
As $f$ is compactly supported, we have $\mu\left(f\right)=lim_{\epsilon\rightarrow0}c_\epsilon\left(f\right)$.
(Here we use that we are admitting ideal simplices: otherwise the support of $f$ would not be compact.)\\
Now using $\lim_{\epsilon\rightarrow 0}c_\epsilon\left(S_\delta\right)=0$ we have $$\mu^\pm\left(f\right)=\lim_{\epsilon\rightarrow 0}c^\pm_\epsilon\left(f\right)\ge\lim_{\epsilon\rightarrow 0}c^\pm_\epsilon\left(SS_n\left(M\right)\right)-c^\pm_\epsilon\left(S_\delta\right)=\lim_{\epsilon\rightarrow0}c^\pm_\epsilon\left(SS_n\left(M\right)\right).$$
But $c^+_\epsilon\left(SS_n\left(M\right)\right)+c_\epsilon^-\left(SS_n\left(M\right)\right)=
\parallel c_\epsilon^+\parallel+\parallel c_\epsilon^-\parallel=\parallel c_\epsilon\parallel\ge\parallel M,\partial M\parallel $
implies that one of $\lim_{\epsilon\rightarrow0}c^\pm_\epsilon\left(SS_n\left(M\right)\right)$ must be at least $\frac{1}{2}\parallel M,\partial M\parallel$, thus positive.\end{pf}

\begin{lem}: Efficient fundamental chains $\mu$ are cycles, i.e.
$\left(\partial\mu\right)^+=\left(\partial\mu\right)^-=0$.\end{lem}
\begin{pf} Denote by $T^i_\epsilon\left(M\right)$
the set of (possibly ideal) $i$-simplices
intersecting $\thick$.                                                               For
all $\delta < \epsilon$, we get by the convexity
of $M_{\left[0,\delta\right]}\subset M_{\left[0,\epsilon\right]}$:
$$ B\subset T^{n-1}_\epsilon\left(M\right)\mbox{ measurable }\Rightarrow
\partial c_{\delta}^\pm\left(B\right)=0.$$
 
When $\partial\mu^\pm$ is a weak-* accumulation point of
a sequence $\partial c_\delta^\pm$,
we conclude
$\partial\mu^\pm\left(B\right)= 0$ for all measurable sets $B$ contained in
{\em some} $T^{n-1}_\epsilon\left(M\right)$ by part (ii) of lemma 2.2, since we may consider them as
subsets of an open set still contained in some slightly larger
$T^{n-1}_\epsilon\left(M\right)$.
 
But clearly, $\cup_{k=1}^\infty T^{n-1}_{\frac{1}{k}}\left(M\right)$ is the
set of all (even ideal) (n-1)-simplices, hence the claim of the lemma.\end{pf}\\     
{\em Remark:} In the case of closed manifolds, this lemma is, of course, an                 immediate consequence of the fact 
that $\parallel\partial\parallel\le n+1$.

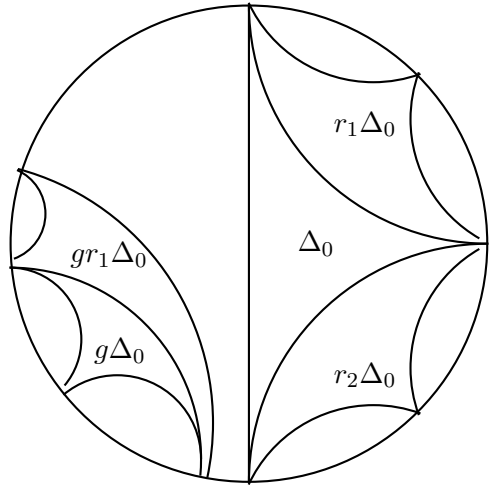
\psset{unit=0.117\hsize}
$$\pspicture(-2,-2)(2,2)
\pscircle(0,0){2}
\psline(0,2)(0,-2)
\psarc(2,2){2}{180}{270}
\psarc(2,-2){2}{90}{180}
\psarc(2.5,1.04){1.15}{160}{240}
\psarc(1.04,2.5){1.15}{205}{290}
\psarc(2.5,-1.04){1.15}{120}{200}
\psarc(1.04,-2.5){1.15}{70}{155}
\psarc(-2,-1.8){1.6}{355}{90}
\psarc(-1.1,-1.8){0.7}{0}{130}
\psarc(-2,-0.8){0.6}{320}{90}
\psarc(-2.5,-1.5){2.2}{348}{75}
\psarc(-2.1,0.25){0.4}{290}{66}

\uput[0](0.3,0){$\Delta_0$}
\uput[0](0.6,1){$r_1\Delta_0$}
\uput[0](0.6,-1.1){$r_2\Delta_0$}
\uput[0](-1.4,-0.9){$g\Delta_0$}
\uput[0](-1.6,-0.1){$gr_1\Delta_0$}
\endpspicture$$       
\subsection{Invariance under ideal reflection group}
Since we have an ordering of the vertices of a simplex $\Delta$, we can speak
of the i-th face of $\Delta$, the codimension 1-face not containing the
i-th vertex.  

\begin{df}: Fix a regular ideal simplex $\Delta_0\subset {\Bbb H}^n$ and, for
$i=0,\ldots,n$, let $r_i$ be the
reflection in the i-th face of $\Delta_0$.
Let $R\subset Isom\left({\Bbb H}^n\right)$
be the subgroup generated by $r_0,\ldots,r_n$ and
let $R^+=R\cap Isom^+\left({\Bbb H}^n\right)$.\end{df}   

We know that $\mu^\pm$
are measure cycles supported on the
set of regular ideal simplices. By 2.1, we may consider $\mu^\pm$
as measures on $\Gamma\backslash Isom\left({\Bbb H}^n\right)$,
after fixing some regular ideal simplex $\Delta_0$ in ${\Bbb H}^n$. 

We will use the convention that $\gamma\in Isom\left({\Bbb H}^n\right)$ 
corresponds to the simplex $\gamma\Delta_0$, i.e., we let $Isom\left({\Bbb H}^n\right)$, and in particular $\Gamma$, act from the {\bf left}. It will be important to note that, after this identification, the {\bf right}-hand action of $R$ corresponds to the following operation
on the set of simplices: $r_i$ maps a simplex to the simplex obtained by reflection in the
i-th face. This is clear from the picture above.

\begin{lem}: For $n\ge3$, efficient fundamental chains are invariant
under the right-hand action
of $R^+$ on $\Gamma\backslash Isom\left({\Bbb H}^n\right)$.\end{lem}
\noindent
Note: If
$\Delta=g\Delta_0$ for some $g\in\Gamma\backslash Isom\left({\Bbb H}^n\right)$,
then the reflection $s_i$ in the
i-th face of $\Delta$ maps $\Delta=g\Delta_0$ to $gr_i\left(\Delta_0\right)$.
In other words, the choice of another reference simplex
changes the identification with $Isom\left({\Bbb H}^n\right)$ by {\bf left}
multiplication with $g\in Isom\left({\Bbb H}^n\right)$, but does
not alter the
{\bf right}-hand action of $R^+$ on $Isom\left({\Bbb H}^n\right)$.
This implies that the truth of lemma 3.9 is independent
of the choice of $\Delta_0$.\\
 
Lemma 3.9
follows from
\begin{lem}: In dimensions $n\ge3$, a signed alternating measure $\mu$ on the set of maximal
volume simplices is
a cycle iff
$r_i^*\left(\mu\right)=-\mu$
for all $i=0,\ldots,n$. \end{lem}
\begin{pf} 
If $n\ge3$, then for any ordered regular ideal (n-1)-simplex $\tau$,
there are exactly two ordered regular ideal n-simplices, $\tau_i^+$ and
$\tau_i^-$,
having $\tau$ as i-th face. (By the way, this is besides lemma 2.6 and its 'corollary' lemma 3.2 the only
point entering the proofs of our theorems which uses
$n\ge3$.)
We fix them such that $\tau_i^+$ is positively oriented.
For a measurable set $B\subset\left\{\mbox{ordered regular ideal
(n-1)-simplices}\right\}$ define

$B_i^+=\left\{\tau_i^+:\tau\in B\right\}$\ \ \ and
$B_i^-=\left\{\tau_i^-:\tau\in B\right\}$.\\
Since $\mu$ is supported on the set of regular ideal n-simplices, we have
that
$$\partial\mu^\pm\left(B\right)=\sum_{k=0}^n\left(-1\right)^k\mu^\pm\left(\partial_k^{-1}\left(B\right)\right)$$
$$=\sum_{k=0}^n\left(-1\right)^k\left(\mu^\pm\left(B_k^+\right)+\mu^\pm\left(B_k^-\right)\right).$$
We may assume that $\mu$ is alternating, in particular $\pi_{ik}^*\mu=\left(-1\right)^{i-k}\mu$, where $\pi_{ik}$ is induced by the affine map realizing the
transposition of the i-th and k-th vertex. 
$\pi_{ik}$ maps $B_i^+$ to $B_k^+$
and $B_i^-$ to $B_k^-$. 
Therefore, for any $i\in\left\{0,\ldots,n\right\}$, we get
$$\partial\mu^\pm\left(B\right)=\sum_{k=0}^n\left(-1\right)^k\left(-1\right)^{i-k}\left(\pi_{ik}^*\mu^\pm\left(B_k^+\right)+\pi_{ik}^*\mu^\pm\left(B_k^-\right)\right)$$
$$=\sum_{k=0}^n\left(-1\right)^i\left(\mu^\pm\left(B_i^+\right)+\mu\left(B_i^-\right)\right)$$
$$=\left(-1\right)^i\left(n+1\right)\left(\mu^\pm\left(B_i^+\right)+\mu^\pm\left(B_i^-\right)\right).$$
In particular $\partial\mu\left(B\right)=0$ holds if and only if
$\mu\left(B_i^+\right)=-\mu\left(B_i^-\right)$ for $i=0,\ldots,n$.
 
The action of $r_i$ maps $B_i^+$ bijectively to
$B_i^-$ and vice versa. This implies the 'if'-part of lemma 3.10.\\
To get the 'only if'-part, we use that $\partial\mu=0$ implies that $r_i^*\mu\left(B_i^\pm\right)=-\mu\left(B_i^\pm\right)$
holds, for any set 
$B\subset \left\{\mbox{ordered regular ideal
(n-1)-simplices}\right\}$.  
Now let $C\subset \left\{\mbox{ordered regular ideal
n-simplices}\right\}$ be an arbitrary set. We divide $C=C^+\cup C^-$, where $C^+=\left\{\sigma\in C: \sigma\mbox{ positively oriented }\right\}$. Consider $B:=\left\{\partial_i\sigma:\sigma\in C^+\right\}$. ($i$ is arbitrary, e.g. $i=0$.) Then we have $B_i^+=C^+$, because for any ordered regular ideal $n-1$-simplex $\partial_i\sigma\in B$, $\sigma$ is the {\em unique} positively oriented ordered regular
ideal $n$-simplex having $\partial_i\sigma$ as its $i$-th face. Thus 
$\mu\left(C^+\right)=-r_i^*\left(C^+\right)$. The same way one gets $\mu\left(C^-\right)=-r_i^*\mu\left(C^-\right)$, 
thus $\mu\left(C\right)=-r_i^*\mu\left(C\right)$. Since $C$ was arbitrary, this proves the 'only if'-part.\end{pf}\\
{\em Remark}: A different
proof of the same fact is given in lemma 2.2.\ 
of \cite{jun}.

\section{Decomposition of efficient fundamental cycles}
 
If $n\ge 4$, then the group generated by reflections in the faces of a regular ideal  
n-simplex in ${\Bbb H}^n$ is dense in $Isom\left({\Bbb H}^n\right)$. We get therefore from lemma 3.9 that
efficient fundamental cycles are invariant under the right-hand action 
of $Isom^+\left({\Bbb H}^n\right)$. This implies that they are a multiple of 
$Haar - r^*Haar$, where $Haar$ is the Haar measure on $Isom^+\left({\Bbb H}^n\right)$. (It is well-known that all invariant measures on a Lie group are multiples of the Haar measure.)

In the following we will discuss the case $n=3$.\\
We wish to recall some facts from the {\bf ergodic theory of unipotent actions}.\\
The Iwasawa decomposition $G=KAN$ of $G=Isom^+\left({\Bbb H}^n\right)$
is as follows:
fix some $v_\infty\in\partial_\infty {\Bbb H}^n$ and some $p\in {\Bbb H}^n$.
Then we may take $K$ to be the group of isometries fixing $p$, $A$ the group
of translations
along the geodesic through $p$ and $v_\infty$,
and $N$ the group of translations
along the horosphere through $p$ and $v_\infty$.\\
We will consider the natural right-hand action of $N$ on $G=KAN$.\\
The next lemma follows from \cite{dan1}. It is nowadays a special case of the Ragunathan conjecture,
which was proved by Ratner.
\begin{lem}: Let G=KAN be the Iwasawa decomposition of
a simple Lie group of ${\Bbb R}$-rank 1, and
$\Gamma\subset G$ a discrete subgroup
of finite covolume. If $\mu$ is a
finite N-invariant ergodic measure on $\Gamma\backslash G$, then
$\mu$ is either a multiple of the Haar measure or it is
supported on a compact $N$-orbit.\end{lem}
The following lemma is a straightforward generalisation of
theorem 4.4. in \cite{dan2}.
\begin{lem}:   Let G=KAN be the Iwasawa decomposition of
a simple Lie group of ${\Bbb R}$-rank 1, and
$\Gamma\subset G$ a discrete subgroup
of finite covolume.
Let $N^\prime\subset N$ be a closed subgroup such that $N/N^\prime$ is compact.
Then any $N^\prime$-invariant ergodic measure on
$\Gamma\backslash G$ is either a multiple of
the Haar measure or is supported on a compact $N$-orbit.\end{lem}
\begin{pf}
We will use several times the following basic fact: If $G_1$ and $G_2$ are subgroups of a group $G$ endowed with
a measure $\mu$, then the left action of $G_1$ on $G/ G_2$ is ergodic if and only if the right action of $G_2$ on 
$G_1\backslash G$ is ergodic. (This is known as Moore-equivalence).

By Moore-equivalence,
ergodic measures for the $N^\prime$-action
on $\Gamma\backslash G$ correspond
to ergodic measures for the action of $\Gamma$ on
$G/N^\prime$.
Consider, therefore, $\mu$ as a measure on $G/N^\prime$,
ergodic with respect to the $\Gamma$-action.
Let $pr:G/N^\prime\rightarrow G/N$ be the projection. Since
$N/N^\prime$ is compact, we have a locally finite
measure $pr_*\mu$ on
$G/N$ which is easily seen to
be ergodic with respect to the $\Gamma$-action. By
lemma 4.1 and Moore-equivalence,
$pr_*\mu$ must either be the Haar measure or correspond to an
$N$-invariant measure
on $\Gamma\backslash G$ which is determined
on a compact orbit $\Gamma\backslash\Gamma gN\subset
\Gamma\backslash G$.

If $pr_*\mu$ = Haar measure, it follows easily that $\mu$
is absolutely continuous with respect to the Haar measure and then
one gets, from ergodicity of the $\Gamma$-action
(theorem 7 in \cite{moo}), that $\mu$ is a
multiple of the Haar measure.

In the second case, $pr_*\mu$ must be
supported on
the $\Gamma$-orbit of some $gN\in G/N$.  Therefore, $\mu$ is
supported on the $\Gamma\times N$-orbit of
$gN^\prime\in G/N^\prime$. By Moore-equivalence
we get a measure supported on the compact $N$-orbit.\end{pf}\\

Lemma 4.1 and 4.2 apply in particular to $G=Isom^+\left({\Bbb H}^n\right)$ with the Iwasawa decomposition described above.

Back to the situation of section 3.2.
Let $v$ be an ideal vertex of the reference simplex
$\Delta_0$. Let $N_v\subset Isom^+\left({\Bbb H}^3\right)$ be the
subgroup of parabolic isometries fixing $v$. We may
consider $N_v$ as the $N$-factor in the Iwasawa decomposition
$Isom^+\left({\Bbb H}^3\right)=K_vA_vN_v$.
(That means we use $v$ and some
arbitrary $p\in {\Bbb H}^3$ to construct the                                    
Iwasawa decomposition.
In the following, we will fix some arbitrary $p\in {\Bbb H}^3$ but consider various
$v\in\partial_\infty {\Bbb H}^3$, therefore the labelling of the
Iwasawa decompositions.)
 
Instead of $R^+$ defined in section 3.2, we consider only the subgroup $T_v^\prime\subset R^+\subset Isom^+\left({\Bbb H}^3\right)$
generated by products of even numbers of
reflections in those faces of $\Delta_0$ which
contain $v$. $\mu$ is, of course,
also invariant under the smaller group $T_v^\prime$.
In \cite{jun} it is shown that $T_v^\prime$ contains a
subgroup $T_v$ which is a
cocompact subgroup of $N_v$ (if $n=3$).
 
The signed measure $\mu$ decomposes as a difference of two (non-negative) measures
$\mu^+$ and $\mu^-$. Both are invariant under the right-hand
action of $T_v$.
From
lemma 2.1, we get that
the probability measures $\overline{\mu}^\pm$,
obtained by rescaling the restrictions of $\mu^\pm$ to
$\Gamma\backslash Isom^+\left({\Bbb H}^3\right)$,
have decomposition maps with respect to the action of
$T_v$,
$$\beta^\pm_v:\Gamma\backslash Isom^+\left({\Bbb H}^3\right)
\rightarrow\mathcal{E}.$$
Here, $\mathcal{E}$ is the set of ergodic $T_v$-invariant measures on
$\Gamma\backslash Isom^+\left({\Bbb H}^3\right)$.
From lemma 4.2, we get that $\mathcal{E}$
consists of $Haar$ (the Haar measure, rescaled to a probability measure) and measures
determined on compact $N_v$-orbits.
The following lemma is well-known.
\begin{lem}: An orbit $gN_v$
is compact in $\Gamma\backslash Isom\left({\Bbb H}^n\right)$ iff
all simplices $gh\Delta_0$ with $h\in N_v$ have its ideal
vertex $g\left(v\right)$ in a parabolic fixed point
of $\Gamma$.\end{lem}
\begin{pf} Parametrise elements of $N_v$ as $u\left(s\right),
s\in {\Bbb R}^{n-1}$ (identifying a stabilized horosphere with
euclidean (n-1)-space). The $N_v$-orbit of $g$ on
$\Gamma\backslash Isom\left({\Bbb H}^n\right)$ is compact if and only if,
for all $s\in {\Bbb R}^{n-1}$, one finds $\gamma\in\Gamma$ and $t\in {\Bbb R}$ such that
$gu\left(ts\right)=\gamma g$.
This $\gamma$ is then conjugated to $u\left(ts\right)$ and, in particular, is
parabolic, i.e., has only one fixed point. The fixed point of $\gamma$ must
be $g\left(v\right)$, since $\gamma g\left(v\right)=gu\left(ts
\right)\left(v\right)=g\left(v\right)$. \\
The other implication is
straightforward.\end{pf}\\

To summarize, we have the following statement:
For any vertex $v$ of the reference simplex $\Delta_0$,
the ergodic
decomposition of the rescaled $\overline{\mu}^\pm$
with respect to the right-hand action of $T_v$
uses the Haar measure and measures determined on
the set of those simplices $g\Delta_0$                                          
which have the vertex $g\left(v\right)$ in a parabolic fixed point of $\Gamma$.
 
\subsection{Manifolds which are not Gieseking-like}
\begin{df}: A 3-manifold is {\bf Gieseking-like} if it has a
hyperbolic structure $M=\Gamma\backslash {\Bbb H}^3$ of finite volume
such that $Q\left(\omega\right)
\cup\left\{\infty\right\}\subset \partial_{\infty}{\Bbb H}^3$
are parabolic fixed points
of $\Gamma$. \end{df}
Here, we have used the upper half space model of ${\Bbb H}^3$,                     and identified
the ideal boundary with ${\Bbb C}\cup\left\{\infty\right\}$. $\omega=\frac{1}{2}+\frac{\sqrt{-3}}{2}$ is 
the 4th vertex of a regular ideal simplex with vertices $0,1,\infty$. The condition is, of course, equivalent to 
the condition that $\Gamma$ is conjugate to a discrete subgroup of $PSL_2Q\left(\omega\right)$ after the 
identification of $Isom^+\left({\Bbb H}^3\right)$ with $PSL_2{\Bbb C}$.
One does not seem to know any example of a Gieseking-like manifold which
is not a finite cover of the Gieseking manifold (communicated to the author
by Alan Reid, see also \cite{lr}). 

\begin{thm}: Let $M$ be a compact manifold of dimension $n\ge3$ such that
$int\left(M\right)$ admits a hyperbolic
metric of finite volume. Assume that $M$ is either of dimension $\ge 4$ or that $M$ is of dimension 3 and is not Gieseking-like.\\
If $\mu$ is an efficient fundamental cycle on $M$, then $\mu=K\left(Haar-r^*Haar\right)$ for some real number $K$.\end{thm}
\begin{pf} By the first remark of chapter 4, we may restrict to dimension 3.
We have to exclude the existence of
a signed measure $\nu$ which is supported on the set of regular ideal simplices with vertices in cusps and which satisfies
$r^*\nu=\pm\nu$ for all $r\in R$. However, the existence of such a nontrivial
signed measure would imply the existence of an $R$-invariant family $\left\{\Delta r:r\in R\right\}$
of simplices with vertices in the cusps of $M=\Gamma\backslash {\Bbb H}^3$. By 2.1, there is 
$g\in Isom\left({\Bbb H}^3\right)$ with $\Delta=g\Delta_0$, where $\Delta_0$ is the ideal simplex with vertices
$0,1,\infty,\omega$ in the upper half-space model. We get that all vertices of the form 
$gv_\infty r$ with $r\in R$ and $v_\infty$ one of $0,1,\infty,\omega$ must be parabolic fixed points of $\Gamma$. 
Note that $\left\{v_\infty r:v_\infty\in\left\{0,1,\omega,\infty\right\},r\in R\right\}=Q\left(\omega\right)\cup\left\{\infty\right\}$. Thus, conjugating $\Gamma$ with $g$ we get a hyperbolic structure with all of $Q\left(\omega\right)\cup\left\{\infty\right\}$
as parabolic fixed points.\end{pf}

\section{Cycles not transversal to geodesic surfaces}
 
In the last section, we classified efficient fundamental cycles on finite-volume hyperbolic manifolds which are not
Gieseking-like. In this chapter, we will see that for arbitrary (possibly Gieseking-like) finite-volume hyperbolic manifolds
we can still obtain information which in chapter 6 will be used to derive glueing inequalities.
\begin{df}: For a hyperbolic manifold $M$ and a two-sided totally geodesic
codimension-1 submanifold $F\subset M$ call \\
- $S^i_{cusp}$ the set of positively oriented ideal
i-simplices with all vertices in parabolic fixed points
of $M$, and\\
- $S^i_F$ the set of (possibly ideal) positively oriented
i-simplices that intersect F transversally.\end{df}
Here, a simplex $\sigma$ is said to intersect $F$ transversally if
it intersects both components of any regular neighborhood of $F$.
\begin{lem}: If $M$ is a hyperbolic manifold and $F$ is a two-sided totally geodesic codimension-1-submanifold, then\\
$S^n_F\cap\left\{ regular\ ideal\ simplices \right\}
\subset\left\{ regular\ ideal\
simplices \right\}$\\ 
has positive Haar measure.\end{lem}
\begin{pf} It is easy to see that
$S^n_F\cap\left\{ regular\ ideal\ simplices \right\}$ is an open, non-empty subset
of $\left\{ regular\ ideal\ simplices \right\}$.\end{pf}
\begin{thm}: Let $M$ be a compact manifold of dimension $n\ge3$ such that
$int\left(M\right)$ admits a hyperbolic
metric of finite volume, and let $F\subset M$ be
a closed totally geodesic codimension-1-submanifold.\\
If $\mu$ is an efficient fundamental cycle (with $\mu^+\mid_{\Gamma\backslash Isom^+\left({\Bbb H}^n\right)}\not=0$),\\
then
$\mu^+\left(S^n_F\right)\not=0$.\end{thm}
\begin{pf} 
Very roughly, the idea is the following: If $\mu^+\left(S_F^n\right)$ vanishes,
then the Haar measure can only give a zero contribution to the ergodic decomposition of $\mu^+$, hence, $\mu^+$ is supported on $S^n_{cusp}$. In particular,
$\mu^+$ vanishes on the set of simplices with boundary faces in $F$, and this will give a contradiction.\\
\\
Rescale $\mu^+\mid_{\Gamma\backslash Isom^+\left({\Bbb H}^n\right)}$ to a probability measure $\overline{\mu}^+$.\\
Assume for some totally geodesic surface $F$ we had
$\overline{\mu}^+\left(S_F^n\right)=\mu^+\left(S^n_F\right)=0$.
 
Let $v$ be a vertex of the reference simplex $\Delta_0$.
Using the ergodic decomposition with respect to the $T_v$-action
on $\Gamma\backslash G=\Gamma\backslash Isom^+\left({\Bbb H}^n\right)$ yields
$$0=\overline{\mu}^+\left(S^n_F\right)=\int_{\Gamma\backslash G}\beta_v\left(g\right)
\left(
S^n_F\right)d\overline{\mu}^+\left(g\right) \ge \int_{g\in\Gamma\backslash G:\beta_v
\left(g\right)=Haar}\beta_v\left(g\right)
\left(S^n_F\right)d\overline{\mu}^+\left(
g\right)$$
$$   =\int_{g\in\Gamma\backslash G:\beta_v
\left(g\right)=Haar}Haar\left(S^n_F\right)d\overline{\mu}^+\left(g\right)
=Haar\left(S^n_F\right)\int_{g\in\Gamma\backslash G:\beta_v
\left(g\right)=Haar}d\overline{\mu}^+\left(g\right)$$
By lemma 5.2,
$ Haar\left(S^n_F\right)\not=0$ and, thus, $$\int_{g\in\Gamma\backslash
G:\beta_v
\left(g\right)=Haar}d\overline{\mu}^+\left(g\right)=0.$$
We will conclude that
$\mu^+$ is supported on
$S^n_{cusp}$ by means of lemma 5.5, which we state separately
because it will be of independent use in chapter 7.                           
 
\begin{df}: Let $\Gamma\subset G=Isom^+\left({\Bbb H}^n\right)$ be a cocompact              discrete subgroup, $v\in\partial_\infty {\Bbb H}^n$, $T_v\subset
Isom^+\left({\Bbb H}^n\right)$
the subgroup defined in chapter 4 and $\beta$ a decomposition map
for the right-hand action of $T_v$, as defined in chapter 4. Define
$$H_v=\left\{g\in \Gamma\backslash G:\beta_v\left(g\right)=
Haar\right\}.$$\end{df}
\begin{lem}: Let $v_0,\ldots,v_n$ be the vertices of a regular ideal simplex
in ${\Bbb H}^n$ and $\overline{\mu}^+$ a probability
measure on $\Gamma\backslash G:= \Gamma\backslash Isom^+\left({\Bbb H}^n\right)$, invariant with respect
to the right-hand action of $R^+$.
If $\overline{\mu}^+\left(H_{v_i}\right)=0$ for
$i=0,\ldots,n$, then $\overline{\mu}^+$ is supported on $S^n_{cusp}$.\end{lem}
\begin{pf}: Let
$$A_i=\left\{g\in\Gamma
\backslash G:gv_i\mbox{ is cusp of }\Gamma\right\}$$ and $$B_i=\left\{g\in
\Gamma\backslash G: \Gamma\backslash \Gamma gN_{v_i}
\mbox{ is compact }\right\}.$$
We have
 
$\Gamma\backslash G - S^n_{cusp}=\Gamma\backslash G - \dura
 = \Gamma\backslash G - \durb=\ver 
$,\\
where the second equality holds by lemma 4.3.\\
 
If $e$ is a $T_{v_i}$-
ergodic measure supported on a compact $N_{v_i}$-orbit,
then
$$e\left(\Gamma\backslash G - B_i\right)=0.$$
Thus (abbreviating $\beta_g:=\beta_{v_i}\left(g\right)$),\\
$$\overline{\mu}^+\left(\Gamma\backslash G - B_i\right)=
\int_{\Gamma\backslash G}\beta_g\left(\Gamma\backslash G - B_i\right)
d\overline{\mu}^+\left(g\right)$$
$$=\int_{H_{v_i}}
\beta_g\left(\Gamma\backslash G - B_i\right)
d\overline{\mu}^+\left(g\right)+\int_{\Gamma\backslash G - H_{v_i}}
\beta_g\left(\Gamma\backslash G - B_i\right)
d\overline{\mu}^+\left(g\right)$$
$$=Haar\left(\Gamma\backslash G - B_i\right)\overline{\mu}^+\left(
H_{v_i}\right)+
\int_{\Gamma\backslash G - H_{v_i}}
\beta_g\left(\Gamma\backslash G - B_i\right)d\overline{\mu}^+\left(g\right)$$
$$=Haar\left(\Gamma\backslash G - B_i\right)\ {\bf 0}\ +
\int_{\Gamma\backslash G - H_{v_i}}{\bf 0}\ d\overline{\mu}^+\left(g\right)=
{\bf 0}$$
and, therefore,
$$\overline{\mu}^+\left(\Gamma\backslash G - S^n_{cusp}\right)=
\overline{\mu}^+\left(\cup_{i=0}^n\Gamma\backslash G - B_i\right)\le\sum_{i=0}^n
\overline{\mu}^+\left(\Gamma\backslash G - B_i\right)=0.$$\end{pf}\\
\\
We are now going to finish the proof of theorem 5.3:\\

We know (from the proof of lemma 3.3) that $\Phi_x\left(c_\epsilon^+\right)\ge
\Phi_x\left(c_\epsilon\right)\ge 1$ for all $x\in\thick$. F is
a closed totally geodesic hypersurface. Therefore $F\subset\thick$
for sufficiently
small $\epsilon$. We conclude $\Phi_x\left(
c_\epsilon^+\right)
\ge1$ for all $x\in F$.                                             
 
For $x\in M$ let $S_x^n$ be the set of straight n-simplices $\Delta$ containing $x$ in their
image. If $x\in F$ is contained in the totally geodesic submanifold $F$, then
$S_x^n$ is the union of the following two sets of simplices:

-simplices in $S_x^n$ which intersect $F$ 
transversally, and

-simplices in $S_x^n$ which have a vertex in $F$.

$\mu^+$ vanishes
on the second set, since it is determined
on $S^n_{cusp}$ and the closed totally geodesic hypersurface $F$ can not have cusps.
Thus, we obtain
$$\mu^+\left(S_x^n\right)=\mu^+\left(S^n_F\cap
S^n_x\right)\le\mu^+\left(S_F^n\right),$$
i.e., it suffices to show that $\mu^+\left(S_x^n\right) > 0$.

For a measure $\mu^+$ on $SS_n\left(M\right)$,
let $\Phi_X\left(\mu^+\right)=\int_{SS_n\left(M\right)}\Phi_x\left(\sigma\right)d\mu^+\left(\sigma\right)$,
where $\Phi_x\left(\sigma\right)$ is given in definition 2.3. Weak-*-convergence implies $\Phi_x\left(\mu^+
\right)=\lim_{\epsilon\rightarrow 0}\Phi_x\left(c_\epsilon^+\right)\ge1$.

On the other hand, $\Phi_x\left(\sigma\right)=0$ if $\sigma\not\in S_x^n$, hence 
$$\Phi_x\left(\mu^+\right)=\int_{S_x^n}\Phi_x\left(\sigma\right)d\mu^+\left(\sigma\right).$$
If $\mu^+\left(S_x^n\right)=0$, then $\Phi_x\left(\mu^+\right)=\int_{S_x^n}\Phi_x\left(\sigma\right)d\mu^+\left(\sigma\right)=0$ (regardless whether $\Phi_x$ is bounded or not), giving a contradiction. Thus $\mu^+\left(S_x^n\right)>0$, implying $\mu^+\left(S_F^n\right)>0$.
\end{pf}\\
Remark: If $\mu^+\mid_{Isom^+\left(H^n\right)}=0$, then $\mu^-\mid_{Isom^+\left(H^n\right)}\not=0$ because of lemma 3.61 and lemma 3.10, and we get 
with an analogous proof $\mu^-\left(S_F^n\right)\not=0$.                                                                          
\section{Acylindrical hyperbolic manifolds}
In this chapter we extend corollary 1 from \cite{jun} to manifolds with cusps.\\
If M is a hyperbolic manifold, define its
convex core to be the smallest closed convex subset $C_M$ of $M$ such that the embedding
$C_M\rightarrow M$ is a homotopy equivalence. $C_M$ is either contained in a geodesic codimension 1 submanifold, or it is 
a codimension 0 submanifold with boundary $\partial C_M$. 
(If $dim\left(M\right)=3$, then $\partial C_M$ is, in general, a pleated surface (see \cite{thu}),
i.e, is almost everywhere totally geodesic and is bent along a family of disjoint geodesics.)
We say that $M$ has totally geodesic
boundary if $C_M$ is homeomorphic to M and $\partial C_M$ is a non-empty totally geodesic submanifold of $M$.
Note that we admit that $C_M$ may have cusps.
If $M$ is an orientable geometrically finite hyperbolic 3-manifold (with $C_M$ not contained in a geodesic codimension 1 submanifold), then $\partial C_M$
is homeomorphic to the union of all non-torus components of the
topological boundary $\partial M$. This applies in particular to any hyperbolic 3-manifold $M$ with totally geodesic boundary $\partial C_M$.
 
Although hyperbolic
structures of infinite
volume are not necessarily rigid, it follows easily from
Mostow's rigidity theorem that on a manifold
of dimension $\ge3$, there can be
at most one hyperbolic metric $g_0$
admitting totally geodesic boundary, up to isometry.
In particular, the volume of the convex core
with respect to the metric $g_0$ is a topological
invariant. Actually, it was shown in \cite{bona} that $g_0$ minimizes the
volume of the convex core among all hyperbolic metrics on M.
\begin{lem}: Let $M$ be a compact 2-manifold with boundary
$\partial M=\partial_0M\cup\partial_1M$, such that $M-\partial_0 M$ admits
an incomplete hyperbolic metric of finite volume
with $\partial_1M$ totally geodesic and the ends corresponding
to $\partial_0M$ complete. Then $$\parallel M,
\partial M\parallel=\frac{1}{V_2}Vol\left(M\right).$$\end{lem}
\begin{pf} It is well-known that any (possibly bounded)
surface of non-positive Euler characteristic
satisfies $\parallel M,\partial M\parallel=-2\chi\left(M\right)$. By the
Gau\ss -Bonnet-formula, this is the same as $\frac{1}{\pi}Vol\left(M\right)=\frac{1}{V_2}Vol\left(M\right)$.\end{pf}
\begin{cor}: Let $n\ge3$ and let $M$ be a compact n-manifold with
boundary $\partial M=\partial_0M\cup\partial_1M$, such                     
that $M-\partial_0 M$ admits
an incomplete hyperbolic
metric of finite volume with $\partial_1M$ totally geodesic
and the ends corresponding to $\partial_0M$ complete. Then, $$\parallel M,\partial M\parallel >
\frac{1}{V_n}Vol\left(M\right).$$\end{cor}
\begin{pf} \\
$\parallel M,\partial M\parallel \ge
\frac{1}{V_n}Vol\left(M\right)$ follows from the familiar argument that
fundamental cycles can be straightened to invoke only simplices of
volume smaller than $V_n$ or, equivalently, from the trivial inequality
$\parallel DM\parallel\le 
2\parallel M,\partial M\parallel$.
 
Suppose we had equality $\parallel M,\partial M\parallel=\frac{1}{V_n}Vol\left(M\right)$. Glue two differently oriented copies of $M$
via $id\mid_{\partial M}$ to get $N=DM$. The incomplete metrics can be
glued along the totally geodesic boundary and, hence, we have
that $N$ is a complete hyperbolic manifold of finite
volume $Vol\left(N\right)=
2 Vol\left(M\right)$. A relative
fundamental cycle for M of norm smaller than $\frac{1}{V_n}Vol\left(M\right)
+\frac{\epsilon}{2}$
fits together with its
reflection to give a relative fundamental
cycle $c_\epsilon$ on $N$ of $l^1$-norm smaller than
$2\frac{1}{V_n}Vol\left(M\right)+\epsilon=
\frac{1}{V_n}Vol\left(N\right)+\epsilon$, but
consisting of simplices which do not intersect
transversally the totally geodesic
surface $\partial M\subset N$, i.e.,
$c^{\pm}_\epsilon\left(S^n_{\partial M}\right)=0$.\\
Recall that a representative $c_\epsilon$ of $\left[N,\partial N\right]$
was used in chapter 3 to get a representative $Str\left(exc\left(r_{\epsilon *}c_\epsilon^\pm\right)\right)$ of 
$\left[int\left(N\right),N_{\left[0,\epsilon\right]}\right]$ with (at most) the same $l^1$-norm. Here $r_\epsilon:\left(N,\partial N\right)\rightarrow
\left(N_{\left[\epsilon,\infty\right]},\partial
N_{\left[\epsilon,\infty\right]}\right)$ was a homeomorphism, $exc$ was the canonical inclusion, and $Str$ means straightening.
$c^{\pm}_\epsilon\left(S^n_{\partial M}\right)=0$ 
implies $Str\left(exc\left(r_{\epsilon *}c_\epsilon^\pm\right)\right)\left(S^n_{\partial M}\right)=0$,
because $r_\epsilon$ can be choosen to be the identity in a neighborhood of the totally geodesic hypersurface $\partial M\subset N$ (which belongs to the thick part if $\epsilon$ is small enough),
and because straightening in $N=DM$ preserves the set of simplices
not intersecting transversally
the totally geodesic surface $\partial M$.
 
By lemma 3.2, we have some accumulation point
$\mu$ of $Str\left(exc\left(r_{\epsilon *}c_\epsilon\right)\right)$ for a sequence of $\epsilon$ tending to
zero.
Similarly to lemma 5.2, it is easy to see that $S^n_{\partial M}$ is open
in $SS_n\left(N\right)$. Hence,
we can apply
part (ii) of lemma 2.2 to get $\mu^+\left(S^n_{\partial M}\right)=0$.
But this contradicts theorem 5.3.\end{pf}
 
\begin{thm}: (a) Let $n\ge3$ and let $M_i, i=1,2$ be compact n-manifolds with
boundaries $\partial M_i=\partial_0M_i\cup\partial_1M_i$, such that
$M_i-\partial_0 M_i$ admit incomplete hyperbolic metrics of finite volume
with $\partial_1 M_i$ totally geodesic and the ends corresponding to
$\partial_0M_i$ complete.
If $\partial_1^\prime M_i\subset\partial_1 M_i$ are non-empty
sets of connected components of $\partial_1M_i$, $f:\partial_1^\prime M_1
\rightarrow\partial_1^\prime M_2 $ is an orientation-reversing isometry, and $M=M_1\cup_f M_2$, then
$$\parallel M,\partial M\parallel < \parallel M_1,\partial M_1\parallel +\parallel
M_2,\partial M_2\parallel.$$
(b) Let $n\ge3$ and let $M_0$ be a compact n-manifold with
boundary $\partial M_0=\partial_0M_0\cup\partial_1M_0$, such that
$M_0-\partial_0M_0$ admits an incomplete hyperbolic metric of
finite volume with $\partial_1M_0$ totally geodesic and the ends corresponding to $\partial_0M_i$ complete.
If $\partial_1^\prime M_0\subset\partial_1
M_0$ is a non-empty
set of connected components of $\partial_1 M_0$, and $f:\partial_1^\prime
M_0\rightarrow\partial_1^\prime M_0$                                                   
is an orientation-reversing isometry of $\partial_1^\prime
M_0$ exchanging
the connected components by pairs, then, letting $M=M_0/f$,
$$\parallel M,\partial M\parallel < \parallel M_0,\partial M_0\parallel.$$\end{thm}
\begin{pf} (a) The incomplete hyperbolic metrics on $M_1$ and
$M_2$ glue together
to give a complete hyperbolic metric on M of volume $Vol\left(M\right)=
Vol\left(M_1\right)+
Vol\left(M_2\right)$. By the Gromov-Thurston theorem, we know that\\
$\parallel M,\partial M\parallel=\frac{1}{V_n}Vol\left(M\right)$ and, by
corollary 6.2, we have $\parallel M_i,\partial M_i\parallel >\frac{1}{V_n}
Vol\left(M_i\right)$. The claim follows.\\
The proof of (b) is similar.\end{pf}\\

If $\partial M_i$ has dimension $\ge 3$, any homeomorphism is homotopic to an isometry,
by Mostow rigidity. We conclude:
\begin{cor}: Let $n\ge4$, and let $M_1,M_2,M_0$  and $\partial^\prime_1M_i$
satisfy all assumptions of
theorem 6.3.\\
If $f:\partial_1^\prime M_1
\rightarrow\partial_1^\prime M_2$ is a homeomorphism, then
 
$\parallel M_1\cup_fM_2,\partial \left(M_1\cup_fM_2\right)\parallel < \parallel
M_1,\partial M_1\parallel +\parallel
M_2,\partial M_2\parallel$.  \\
If $f:\partial_1^\prime M_0\rightarrow\partial_1^\prime M_0$
is an orientation-reversing homeomorphism of $\partial_1^\prime M_0$ exchanging
the boundary components by pairs, then
 
$\parallel M_0/f,\partial \left(M_0/f\right)\parallel <
\parallel M_0,\partial M_0\parallel$.\end{cor}

\section{Branching of laminations}

We are going to extend results of \cite{cal} to manifolds with cusps (which are not Gieseking-like).\\ 
In this chapter, we always consider foliations/laminations of codimension 1.\\
For more background on the Gromov norm of foliations (and foliations in
general), we refer to \cite{cal}.
\begin{df}: Let $M$ be a manifold, possibly with boundary, and
$\F $ a lamination of M. Define
$$\parallel M,\partial M\parallel_{\F}:=
inf\left\{\sum_{i=1}^r\mid
a_i\mid:\sum_{i=1}^ra_i\sigma_i\in\left[M,\partial M\right],
\sigma_i\mbox{ transverse to }\F\right\}.$$\end{df}
Here, a simplex $\sigma$ is said to be transverse
to the lamination $\F$, if the induced lamination
$\F\mid_\sigma$ is topologically conjugate to the
subset of a foliation of $\sigma$
by level sets of an affine map $f:\sigma\rightarrow R$.

A typical example for non-transversality
of a tetrahedron $\Delta$ to a lamination $\F$
is the following: let $e_1,e_2,e_3$ be the three edges of a
face $\tau\subset\Delta$.
If $\F\mid_{\tau}$ contains three lines which connect respectively
$e_1$ to $e_2$, $e_2$ to $e_3$ and $e_3$ to $e_1$, then $\Delta$
can npt be transverse to $\F$.
 
{\em Remark:} If $\F$ is not transverse to $\partial M$ nor
contains $\partial M$ as a leaf, then
$\parallel M,\partial M\parallel_{\F} =\infty$. Otherwise the
foliated Gromov norm is finite. In the following, we will always assume that
either $\F$ is transverse to $\partial M$ or that $\partial M$ is a leaf of $\F$.\\
\subsection{Laminations with no or one-sided branching}
Recall that a codimension 1 {\bf lamination} $\mathcal{F}$ of an n-manifold $M$
is a decomposition of a closed subset $\lambda\subset M$ into codimension 1 submanifolds (leaves) so that $M$ is covered by charts of the
form $I^{n-1}\times I$, the intersection of a leaf with a chart being of the form $I^{n-1}\times\left\{*\right\}$. A lamination $\mathcal{F}$ of a 3-manifold $M$ with image $\lambda\subset M$
is called {\bf essential} if no leaf is a sphere or a
torus bounding a solid torus, $\overline{M-\lambda}$ is irreducible, and $\partial\overline{M-\lambda}$ is 
incompressible and end-incompressible in $\overline{M-\lambda}$, where the closure of $M-\lambda$ is taken w.r.t.\ any metric (see \cite{goe}, ch.1). E.g., if $M$ is a 3-manifold and $\mathcal{F}$ a foliation without Reeb components, this is an essential lamination.\\
To motivate the following results about the relation between foliated Gromov norm and branching of laminations,
we recall a notion from \cite{goe}. An {\bf order tree} is a set $T$ together with a collection $S$ of linearly ordered segments $\sigma$,
each having distinct least and greatest elements, $e\left(\sigma\right)$ and $i\left(\sigma\right)$, respectively, satisfying the following conditions:\\
- if $\sigma\in S$, then $-\sigma\in S$ (the set $\sigma$ with reversed order),\\
- any closed subset of $\sigma\in S$ belongs to $S$,\\
- any $v,w\in T$ can be joined by some $\sigma_1,\ldots,\sigma_k\in S$, i.e., $v=i\left(\sigma_1\right), e\left(\sigma_j\right)=i\left(\sigma_{j+1}\right)$ and $w=e\left(\sigma_k\right)$,\\
- any cyclic word $\sigma_0,\ldots,\sigma_{k-1}$ with $e\left(\sigma_j\right)=i\left(\sigma_{j+1}\right)$ for all $j$ and $e\left(\sigma_{k-1}\right)=i\left(\sigma_0\right)$ becomes trivial after subdividing the $\sigma_i$'s and performing trivial cancellations,\\
- if $\sigma_1\cap\sigma_2=\left\{i\left(\sigma_2\right)\right\}=\left\{e\left(\sigma_1\right)\right\}$, then $\sigma_1\cup\sigma_2\in S$.\\
To a codimension 1 lamination $\mathcal{F}$ of $M$, one considers the pull-back lamination $\widetilde{\mathcal{F}}$ of $\widetilde{M}$
with image $\lambda\subset \tilde{M}$, and
one constructs $\left(T,S\right)$ as follows:\\
elements of $T$ are either leaves of $\widetilde{\mathcal{F}}$ not contained in $\overline{\widetilde{M}-\lambda}$ or components of $\widetilde{M}-\lambda$.
For each directed arc $\alpha$ intersecting leaves of $\widetilde{\mathcal{F}}$ 
transversally, having non-empty intersection with at least two leaves and being effective (no two points on the intersection with a leaf can be cancelled by an obvious homotopy), let $\sigma$ be the 
set of elements (of $T$) intersected by $\alpha$, with the inherited linear order.\\
According to \cite{goe}, prop.6.10, $\left(T,S\right)$ is an order tree if $\mathcal{F}$ is an essential lamination.

$\mathcal{F}$ is then called ${\Bbb R}$-covered, one-sided branched, or two-sided
branched according to whether the leaf space of $\widetilde{\F}$, considered as an order tree,
is ${\Bbb R}$, branched in one direction,
or branched in both directions. 
 
\begin{lem}: If
$\F$ is an ${\Bbb R}$-covered or one-sided branched
essential lamination on a 3-manifold $M$ such that ${\mathcal{F}}\mid_{\partial M}$ is ${\Bbb R}$-covered, then
$$\parallel M,\partial M\parallel=\parallel M,\partial
M\parallel_{\F}.$$
\end{lem}                                                           
\begin{pf} This is shown in theorems 2.2.10 and 2.5.9 of \cite{cal},
assuming that $M$ is closed. However, the proof
works also for manifolds with boundary.
 
Indeed, since
$\partial M$ is either
transverse to $\F$ or is a leaf of $\F$, the
straightening defined in lemma 2.2.8 of \cite{cal}, for chains with vertices on comparable leaves, preserves
$C_*\left(\partial M\right)$. This implies, in particular, the
claim for ${\Bbb R}$-covered foliations.
In the case of one-sided branching (say in positive direction), the argument in 2.5.9 of
\cite{cal} was then
to isotope a chosen lift of the finite
singular chain in $\tilde{M}$ in the negative direction until its
vertices are on comparable leaves. (This has to be done $\pi_1M$-equivariantly
in the sense that the projection to $M$ stays a relative cycle.)
If $\partial M$ is a leaf of
$\F$, then one can leave all vertices on $\partial M$ fixed and
only isotope the other vertices. If $\partial M$ is transversal to
$\F$, the isotopy can clearly be performed in such a way that
vertices on $\partial M$ (which already are on comparable leaves since ${\mathcal{F}}\mid_{\partial M}$ is ${\Bbb R}$-covered)
are isotoped inside $\partial M$.
 
Hence, in any case,
the straightening maps $C_*\left(\partial M\right)$ to
$C_*\left(\partial M\right)$ and, by the five lemma, it induces
the identity map in relative homology. Thus, it maps
relative fundamental cycles to
relative fundamental cycles transversal to $\F$, not
increasing the $l^1$-norm. \end{pf}
\subsection{Asymptotically separated laminations}
\begin{df}: Let $int\left(M^n\right)$ be hyperbolic and let $\F$ be
a
lamination of $M$.
Let $\widetilde{\F}\mid_{int\left(\widetilde{M}\right)}$ be the pull-back
lamination of ${\Bbb H}^n$. $\F$
is
called asymptotically separated if, for some leaf $F\in\widetilde{\F}$,
there are two geodesic (n-1)-planes
on distinct sides of $F$.\end{df}
We include a proof of the following lemma, implicit in \cite{cal}, for
lack of an explicit reference and
because it might help to understand the idea behind
theorem 7.5.
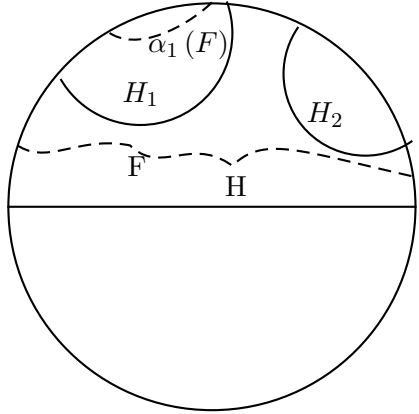
\psset{unit=0.1\hsize}
$$\pspicture(-2,-2)(2,2)
\pscircle(0,0){2}
\psarc(1.5,1.3){0.8}{145}{305}
\psarc(-0.7,1.7){0.9}{210}{20}
\psline(-2,0)(2,0)
\psbezier[linestyle=dashed](-1.9,0.6)(-1.6,0.4)(-1.3,0.7)(-0.8,0.6)(-0.5,0.3)(-0.2,0.7)(0.2,0.4)(0.5,0.7)(0.85,0.55)(1.95,0.3)
\psbezier[linestyle=dashed](-1,1.7)(-0.7,1.5)
(-0.2,1.8)(0,2)
\uput[0](0,0.2){H}
\uput[0](-0.95,0.4){F}
\uput[0](-1,1.1){$H_1$}
\uput[0](0.8,0.9){$H_2$}                                           
\uput[0](-0.75,1.58){$\alpha_1\left(F\right)$}
\endpspicture$$
\begin{lem}: If $\F$ is an asymptotically separated lamination of
a finite-volume hyperbolic manifold $M=\Gamma\backslash {\Bbb H}^n$, then $\F$ is two-sided
branched.\end{lem}
\begin{pf}
Let $F$ be a leaf of
$\widetilde{\F}$ such that
there exist geodesic (n-1)-planes on distinct sides of $F$. These two planes cut out two 
half-spaces $U_1$ and $U_2$ on distinct sides of $F\subset{\Bbb H}^n$.
Let $H$ be the complement of $U_1$ and let
$H_1$ and $H_2$ be disjoint half-spaces in $U_2$. Note that $F\subset H$.
 
If $\Gamma\backslash{\Bbb H}^n$ has finite
volume, then it is well-known that the $\Gamma$-orbits on
the space of pairs of {\em distinct} points in $\partial_\infty {\Bbb H}^n$
are dense.
 
In particular, fixing some arbitrary $\gamma\in\Gamma$ with fixed points
$p_1,p_2$, one finds conjugates of $\gamma$ in $\Gamma$, such that
their fixed points come arbitrarily close to two given points $q_1\not=q_2$
in $\partial_\infty {\Bbb H}^n$. (Namely, conjugate with elements of $\Gamma$
which map $p_1$ close to $q_1$ and $p_2$ close to $q_2$.)
 
It follows
that, in a finite-covolume subgroup $\Gamma\subset
Isom^+\left({\Bbb H}^n\right)$, to any given disk $D\subset\partial_\infty {\Bbb H}^n$, one finds loxodromic isometries with
both fixed points in this disk. Let $\alpha_1$ resp.\ $\alpha_2$ be
such loxodromic isometries with both fixed points in $\partial_\infty H_1$ resp.
both fixed points in $\partial_\infty H_2$.
Loxodromic isometries map any set in the complement of
a neighborhood of the repelling fixed point, after sufficiently many iterations, inside any neighborhood of the attracting fixed point.
Hence, replacing $\alpha_1$ and $\alpha_2$ by sufficiently large powers, we
get that $\alpha_1\left(F\right)\subset H_1$ and $\alpha_2\left(F\right)
\subset H_2$.
 
Since $\widetilde{\F}$ is $\Gamma$-invariant, we have found incomparable leaves $\alpha_1\left(F\right)$ and $\alpha_2\left(F\right)$ above $F$ and, by analogous arguments,
we also get incomparable leaves below $F$. \end{pf}
\\
{\em Remark}: A conjecture of Fenley would
imply
that a foliation of a finite-volume hyperbolic 3-manifold $int\left(M\right)$ is two-sided branched if and only if it is
asymptotically separated, see the discussion in chapter 2.5.\ of \cite{cal}. Namely,
Calegari proves that a two-sided branched foliation (on an arbitrary
hyperbolic manifold) either is
asymptotically separated or the leaves have as limit sets all of
$\partial_\infty {\Bbb H}^3$. On the other hand, Fenley conjectures that for
foliations of finite-volume hyperbolic manifolds (which are transversal to
the boundary $\partial M$), the limit set of a leaf can be all of $\partial_\infty {\Bbb H}^3$
only if $\F$ is ${\Bbb R}$-covered.\\
\\
\\
The following theorem 7.5 extends theorem 2.4.5 in \cite{cal} to the cusped case.\\
\begin{thm}: If the interior of M is a hyperbolic n-manifold of finite
volume
which is not Gieseking-like, $n\ge 3$, and
if $\F$ is an asymptotically
separated lamination, then
$$\parallel M,\partial M\parallel<\parallel M,\partial
M\parallel_{\F}.$$\end{thm}
\begin{pf}
\\
We want to give an outline of the proof. We will show that there exist three 
half-spaces $D_0, D_1, D_2$ such that the following holds: whenever a straight simplex has at least one vertex in each of $D_0, D_1, D_2$, it can not be transverse
to $\F$. 
Assuming $\parallel M,\partial M\parallel_{\F}=\parallel M,\partial M\parallel$, we would have
an efficient fundamental cycle $\mu$ which actually comes from a sequence of
fundamental cycles transverse to $\F$. If $M$ is closed, one gets easily that
$\mu^\pm$ have to vanish on the set of those ideal simplices with at least one vertex 
in each of $\partial_\infty D_0, \partial_\infty D_1, \partial_\infty D_2$. If 
$M$ has cusps, we still get the slightly weaker statement that $\mu^\pm$ have to vanish on the set of those ideal simplices with at least one vertex
in each of $\partial_\infty D_0 - P, \partial_\infty D_1 - P, \partial_\infty D_2 - P$, where $P$ is the set of parabolic fixed points of $\Gamma$. We can then use 
our knowledge of $\mu$ to derive a contradiction.\\

Let $F$ be a leaf which has the property in
the definition of "asymptotically
separated", i.e., there are planes, and hence half-spaces $U_1$ and $U_2$, in disjoint
components of ${\Bbb H}^n - F$. We choose in $U_2$ two smaller disjoint half-spaces
$H_1$ and $H_2$.
Like in the proof of
lemma 7.4, one finds loxodromic isometries $\alpha_1\in\Gamma$ with both fixed points
in $H_1$ and $\alpha_2\in\Gamma$ with both fixed points in $H_2$.
Replacing, if necessary, $\alpha_1$ and $\alpha_2$ by sufficiently large powers,
we arrange that
$\alpha_1\left(U_1\right)\subset H_1$ and $\alpha_2\left(U_1\right)\subset
H_2$,
and that $F,\alpha_1\left(F\right),\alpha_2\left(F\right)$ are disjoint.
Letting $D_0=U_2, D_1=\alpha_1\left(U_1\right)$, and,
$D_2=\alpha_2\left(U_1\right)$,
the remark after definition 7.1 tells us that there is no
tetrahedron transverse to $\tilde{\F}$
with one vertex in each of $D_0,D_1$ and $D_2$.

For the convenience of the reader, we first explain the proof for {\bf closed manifolds}.
Assume that we have straight fundamental cycles $c_\epsilon$,
transverse to $\F$, with 
$\parallel c_\epsilon\parallel
<\parallel M\parallel+\frac{\epsilon}{V_n}$, and that $\mu$ is the
weak-*-limit of $c_\epsilon$. Denoting by $S_{D_0,D_1,D_2}$ the open set of straight
(possibly ideal) simplices with one vertex in each of $D_0,D_1$ and $D_2$, we have just seen that transversality to
$\F$ implies $c_\epsilon^\pm\left(S_{D_0,D_1,D_2}\right)=0$. This implies $\mu^\pm
\left(S_{D_0,D_1,D_2}\right)=0$, contradicting the fact that $\mu^+$ is the Haar measure.
(A similar argument was given by Calegari in 2.4.5 of \cite{cal}.)
 
Now we are going to consider {\bf hyperbolic manifolds of finite volume}.
Let $P\subset\partial_\infty {\Bbb H}^n$ be the parabolic fixed points of
$\Gamma$, and $H_\epsilon=p^{-1}\left(M_{\left[0,\epsilon\right]}\right)
\subset {\Bbb H}^n$ the preimage of the $\epsilon$-thin part. It is the union of horoballs centered at the points of $P$.
For $\delta$ sufficiently small, $D_0 - \overline{H_\delta},
D_1 - \overline{H_\delta}$ and
$D_2 - \overline{H_\delta}$
are nonempty. Fix such a $\delta$.
Let $$S_{D_0,D_1,D_2}=\left\{\mbox{ 
simplices having
vertices } v_0\in D_0 - \overline{H_\delta}, v_1\in D_1 - \overline{H_\delta},
v_2\in D_2 - \overline{H_\delta}\right\},$$ where we admit ideal simplices.\\
We have seen that simplices in
$S_{D_0,D_1,D_2}$ are not transversal to $\tilde{\mathcal{F}}$.                            Moreover, we define $$Str\left(S_{D_0,D_1,D_2}\right):=\left\{Str\left(\sigma\right):\sigma\in S_{D_0,D_1,D_2}\right\}$$

and $${\Bbb U}:=\left\{\mbox{ pos.\ or.\ regular
ideal simplices }\left(v_0,\ldots,v_n\right):
v_i\in \partial_\infty D_i - P\mbox{ for }i=0,1,2\right\}.$$  
 
Now suppose we had the equality $\parallel M,\partial M\parallel=
\parallel M,\partial M\parallel_{\mathcal{F}}$. We will stick to the
notations of chapter 3. Take some
transverse relative fundamental
cycle $c_\epsilon$
of norm smaller than $\parallel M,\partial M\parallel + \frac{\epsilon}{V_n}$
and
make it, via the homeomorphism $r_\epsilon:\left(M,\partial M\right)\rightarrow\left(\thick,\partial\thick\right)$, to a relative
fundamental cycle $r_{\epsilon *}\left(c_\epsilon\right)$
of the $\epsilon$-thick part, which is transverse to
the foliation $r_\epsilon\left(\mathcal{F}\right)$.
We may arrange $r_\epsilon$ to be the identity on the $\epsilon^\prime$-thick
part for $\epsilon^\prime$ slightly larger than $\epsilon$. Then, the lift of
$r_{\epsilon *}\left(c_\epsilon\right)$ to ${\Bbb H}^n$ is
transverse to $\widetilde{\mathcal{F}}$ outside
$H_{\epsilon^\prime}$. By choosing $\epsilon$ and $\epsilon^\prime$
sufficiently small,
one may make this exceptional set $H_{\epsilon^\prime}$
as small as one wishes.\\
Decompose $S_{D_0,D_1,D_2}$ as a
countable union $S_{D_0,D_1,D_2}=\cup_{i=1}^\infty V_i$, where $V_i\subset S_{D_0,D_1,D_2}$ is the open subset of
(possibly ideal) positively oriented
simplices $\sigma\in S_{D_0,D_1,D_2}$ satisfying
$\sigma\cap H_{\frac{1}{i}}=\emptyset$. (The union is all of
$S_{D_0,D_1,D_2}$ because any ideal
or non-ideal simplex
with vertices outside $H_\delta$ must remain outside some $H_{\frac{1}{i}}$
for sufficiently large $i$.)
Let $W_i=Str\left(V_i\right)=\left\{ str\left(\sigma\right):\sigma\in V_i\right\}$.
For $\epsilon$ sufficiently small (such that we can choose $\epsilon^\prime<\frac{1}{i}$),
we have $r_{\epsilon *}\left(c_\epsilon^\pm\right)
\left(V_i\right)=0$, since $r_{\epsilon *}\left(c_\epsilon\right)$ is transverse
to $\mathcal{F}$ outside $H_\frac{1}{i}$ and $V_i$ consists of simplices
which do not intersect $H_\frac{1}{i}$ and which
are not transverse to $\mathcal{F}$.
As a consequence, $Str\left(exc\left(r_{\epsilon *}\left(c_\epsilon\right)\right)\right)\left(W_i\right)=0$.
If $\mu$ is the weak-* limit of the sequence $Str\left(exc\left(r_{\epsilon *}\left(c_\epsilon\right)\right)\right)$ with
$\epsilon\rightarrow 0$, we get
$\mu^\pm\left(W_i\right)=0$ by openness of $W_i$ and part (ii) of lemma 2.2.
 
$W=Str\left(S_{D_0,D_1,D_2}\right)=\left\{Str\left(\sigma\right):\sigma\in S_{D_0,D_1,D_2}\right\}$
is a countable increasing union $W=\cup_{i=1}^\infty W_i$.
Hence $\mu^\pm\left(W\right)=0$.
${\Bbb U}\subset W$ implies $$\mu^\pm\left({\Bbb U}\right)=0.$$
 
On the other hand, ${\Bbb U}$ has nontrivial Haar measure. Indeed, $Isom^+\left({\Bbb H}^n\right)$ corresponds to ordered n-tuples of points in $\partial_\infty {\Bbb H}^n$, because
any
such ordered n-tuple is the set of first $n$ vertices for some unique positively oriented
ordered
regular ideal simplex. Hence, the set of positive regular ideal simplices,
with $v_i\in\partial_\infty D_i$ for $i=0,1,2$, corresponds to an open set
of positive Haar measure in $Isom^+\left({\Bbb H}^n\right)$. Clearly,
a discrete subgroup of
$Isom^+\left({\Bbb H}^n\right)$ has a countable number of parabolic fixed points.
Thus, ${\Bbb U}$ has positive Haar measure.
 
Recall the notation from                    
chapter 4: $v\in\partial_\infty {\Bbb H}^n$ is
an arbitrary vertex of the reference simplex $\Delta_0$
and $\beta_v\left(g\right)$ is the ergodic component of $g\in\Gamma
\backslash Isom^+\left({\Bbb H}^n\right)$ with respect to the $T_v$-action.
We define
$$H_v=\left\{g\in\Gamma\backslash G:\beta_v\left(g\right)=Haar\right\}.$$
 
$Haar\left({\Bbb U}\right)\not=0$ implies $\mu^\pm\left(H_v\right)=0$.
Indeed, from lemma 4.2 and lemma 4.3 we know that
the complement of $H_v$ in the set of regular ideal
simplices is the set of simplices $g\Delta_0$
with the vertex $g\left(v\right)$ in a
parabolic fixed point of $\Gamma$. $\Gamma$ has a countable
number of parabolic fixed points and, therefore, this complement is
a set of trivial Haar measure. Thus,
$$Haar\left({\Bbb U}\cap H_v\right)=Haar\left({\Bbb U}\right)>0$$
and we
apply the ergodic decomposition from 2.2 to get
$$0=\mu^\pm\left({\Bbb U}\cap H_v\right)=
Haar\left({\Bbb U}\cap H_v\right)\mu^\pm\left(H_v\right)$$
which implies $$\mu^\pm\left(H_v\right)=0.$$
 
This discussion applies to all vertices $v_i$ of the reference
simplex $\Delta_0$.
By lemma 5.5, we can conclude that $\mu^\pm$ are determined on $S^n_{cusp}$.
 
In
particular, since $\mu\not=0$, there necessarily {\em are} regular simplices with all vertices
in parabolic fixed points. By lemma 3.10, $\mu$ is invariant up to sign
under the right-hand
action of the regular ideal reflection group $R$ defined in section 3.2.
Hence, there must even be
an $R$-invariant family of regular ideal simplices with vertices in
parabolic fixed points. This is only possible in dimension 3 and, after conjugating with an isometry,
$Q\left(\omega\right)\cup\left\{\infty\right\}$ must be parabolic
fixed points of
$\Gamma$.\end{pf}\\
 
A surface $F$ in a 3-manifold $M$ is called a virtual fiber if there is some
finite cover $p:\overline{M}\rightarrow M$ and some fibration $\overline{F}
\rightarrow\overline{M}\rightarrow S^1$ with $\overline{F}$ isotopic to
$p^{-1}\left(F\right)$. \\
A theorem of Thurston and Bonahon asserts that a
properly embedded compact $\pi_1$-injective surface in a finite-volume
hyperbolic 3-manifold is either quasigeodesic or a virtual fiber.
Since quasigeodesic surfaces are asymptotically separated, one gets analogously to \cite{cal}, theorem 4.1.4:
\begin{cor}: If the interior of $M$ is a
hyperbolic 3-manifold of finite volume which is
not Gieseking-like and $F\subset M$ is a properly embedded compact
$\pi_1$-injective
surface, then
$F$ is a virtual
fiber if and only
if $\parallel M,\partial M\parallel_{\F}
=\parallel M,\partial M\parallel$. \end{cor}

\end{document}